\documentclass[hypertex]{amsart}
\usepackage{amsmath}
\usepackage{amsthm}
\usepackage{amsfonts}
\usepackage{amssymb}

\usepackage{hyperref}

\newlength{\abstand}
\catcode`\ä=\active
\def ä{\"a}
\catcode`\ö=\active
\def ö{\"o}
\catcode`\ü=\active
\def ü{\"u}
\catcode`\Ä=\active
\def Ä{\"A}
\catcode`\Ö=\active
\def Ö{\"O}
\catcode`\Ü=\active
\def Ü{\"U}
\catcode`\ß=\active
\def ß{\ss{}}
\catcode`\é=\active
\def é{\'{e}}
\catcode`\É=\active
\def É{\'{E}}
\catcode`\ê=\active
\def ê{\^{e}}
\catcode`\Ê=\active
\def Ê{\^{E}}

\def\R{{\mathbb R}}

\def\N{{\mathbb N}}

\def\Q{{\mathbb Q}}
\def\Z{{\mathbb Z}}

\def\Ob{\mathrm{Ob}}
\def\cA{\mathcal A}
\def\cB{\mathcal B}
\def\cC{\mathcal C}
\def\cD{\mathcal D}
\def\cF{\mathcal F}
\def\cI{\mathcal I}
\def\cJ{\mathcal J}
\def\cK{\mathcal K}

\def\cP{\mathcal P}
\def\cR{\mathcal R}
\def\cU{\mathcal U}
\def\cV{\mathcal V}

\newcommand\str[1]{{\mbox{}^*#1}}
\newcommand\Mor[3]{\mathrm{Mor}_{#1}(#2,#3)}
\newcommand\id[1]{\mathrm{id}_{#1}}
\newcommand\cov[2]{\mathrm{Funct}\,(#1,#2)}
\newcommand\con[2]{\mathrm{Funct}^\circ(#1,#2)}
\newcommand\Ens[1]{(#1-\underline{\mathrm{Ens}})}
\def\Ensfin{{\mathrm{Ens}}^{\mathrm{fin}}}
\newcommand\Ab[1]{(#1-\underline{\mathrm{Ab}})}
\newcommand\fin[1]{#1^{\text{fin}}}
\def\finG{\fin{\tilde{G}}}
\def\strfin{$\str{\text{finite}}$}
\newcommand\RF[2]{\mathrm{R}^{#2}{#1}}
\newcommand\LF[2]{\mathrm{L}^{#2}{#1}}
\newcommand\Koh[1]{\mathrm{H}\left(#1\right)}

\swapnumbers
\theoremstyle{definition}
\newtheorem{defi}{Definition}[section]
\newtheorem{bsp}[defi]{Example}

\newtheorem{lemma}[defi]{Lemma}

\newtheorem{bem}[defi]{Remark}
\newtheorem{satz}[defi]{Proposition}
\newtheorem{thm}[defi]{Theorem}
\newtheorem{cor}[defi]{Corollary}
\newtheorem*{nonumberlemma}{Lemma}
\newtheorem*{nonumbersatz}{Proposition}

\input xy
\xyoption{all}

\begin{document}

\title{Enlargements of Categories}
\author{Lars Brünjes, Christian Serpé}
\date{\today}
\thanks{The first author is supported by the European Union.}
\thanks{This work has been written under hospitality of the Department of Pure
  Mathematics and Mathematical Statistics of the University of Cambridge and
  the "Sonderforschungsbereich geometrische Strukturen in der Mathematik"
  of the University of Münster}
\address{
  Lars Brünjes \\
  NWF I -- Mathematik \\
  Universität Regensburg \\
  93040 Regensburg \\
  Germany}
\email{lars.bruenjes@mathematik.uni-regensburg.de}
\address{
  Christian Serp\'{e} \\
  Westfälische Wilhelms-Universität Münster \\
  Mathematisches Institut \\
  Sonderforschungsbereich 478 ``Geometrische Strukturen in der Mathematik'' \\
  Hittorfstr. 27 \\
  D-48149 Münster \\
  Germany}
\email{serpe@math.uni-muenster.de}
\date{\today}

\begin{abstract}
In order to apply nonstandard methods to modern algebraic geometry, as a first step in this paper we study the
applications of nonstandard constructions to category theory. It turns out that many categorial properties are
well behaved under enlargements.

\end{abstract}

\maketitle


\tableofcontents

\def\A{\mathbb{A}}

\section{Introduction}

\vspace{\abstand}

After Abraham Robinson's pioneering work on nonstandard analysis in the 1960s (see \cite{nsanalysis}),
nonstandard methods have been applied successfully in a wide area of mathematics,
especially in stochastics, topology, functional analysis, mathematical physics and
mathematical economy.

Nonstandard proofs of classical results are often conceptually more satisfying than the standard
proofs; for example, existence of Haar measures on
locally compact abelian groups can be proved elegantly by nonstandard
methods (see \cite{nshaar}, \cite{nsharmonicanalysis} and \cite{nshaarint}).

Robinson himself also proved the usefulness of his methods for classical problems in
algebraic geometry and number theory,
namely the distribution of rational points on curves over number fields
(see \cite{robnspoints} and \cite{robsiegelmahler}),
infinite Galois theory, class field theory of infinite extensions
(see \cite{robarithmetic} and \cite{robnsnt})
and the theory of Dedekind domains (see \cite{robnsdedekind}).

There are plenty of reasons which could make nonstandard mathematics attractive for arithmetic geometry:
If $\str{\Z}$ is an enlargement of $\Z$, then it contains infinite numbers and in particular infinite prime numbers.
If $P$ is such an infinite prime number, then the residue ring $\str{\Z}/P\str{\Z}$ is a field
that behaves in many respects like a finite field, even though it is of characteristic zero and can --- for
suitable values of $P$ --- even contain an algebraic closure of $\Q$.

If $\str{\Q}$ denotes the field of fractions of $\str{\Z}$, then all $p$-adic fields $\Q_p$ (for standard primes $p$),
as well as the field $\R$ and the ring of adeles $\A_{\Q}$, are all
natural subquotients of $\str{\Q}$
(and the obvious analogues remain true if we replace $\Q$ with a number field $K$ and $\str{\Q}$ with $\str{K}$).
This seems to suggest that working with $\str{\Q}$ could lead to new insights into problems concerning local-global
principles and adelic questions.

Limits, especially limits of ``finite objects",
play an important role in arithmetic geometry: Galois groups are limits of finite groups,
the algebraic closure of a field $k$ is the limit over finite extensions of $k$, Galois cohomology is the limit of
the cohomology of the finite quotients of the Galois group,
$l$-adic cohomology is the limit of \'{e}tale cohomology with coefficients
in the finite sheaves $\Z/l^n\Z$, and so on.
This is another aspect that makes nonstandard methods attractive for arithmetic geometry,
because those methods often make it possible to replace such a limit with a $\str{\text{finite}}$ object
which behaves just like a finite object.
This is exactly the point of view Robinson takes in his approach to infinite Galois theory which allows him
to deduce infinite Galois theory from finite Galois theory. For more motivation for using nonstandard
methods in arithmetic algebraic geometry we refer to \cite{Fes}.

Modern arithmetic geometry makes heavy use of categorical and homological methods,
and this unfortunately creates a slight problem with applying nonstandard methods directly:
The nonstandard constructions are set-theoretical in nature and in the past have mostly been applied to sets
(with structure) like $\R$ or $\Q$ or topological spaces and not to categories and (homological) functors.
But in order to be able to talk about ``$\str{\text{varieties}}$ over $\Q$" or ´"\'{e}tale cohomology
with coefficients in $\str{\Z}/P\str{\Z}$ for an infinite prime $P$",
it seems necessary to apply the nonstandard construction to categories like the category of varieties over $\Q$
or the category of \'{e}tale sheaves over a base scheme $S$.

Exacly for this reason, in this paper, we study "enlargements of categories", i.e. the application of the nonstandard construction
to (small) categories and functors, to create the foundations needed to make use of the advantages of nonstandard
methods described above in this abstract setup.

As a first application in a second paper,
we will enlarge the fibred category of \'{e}tale sheaves over schemes
and then be able to define \'{e}tale cohomology with coefficients in $\str{\Z}/P\str{\Z}$ for an infinite prime $P$
or with coefficients in $\str{\Z}/l^h\str{\Z}$ for a finite prime $l$ and an infinite natural number $h$.
In the case of smooth and proper varieties over an algebraically closed field,
the first choice of coefficients then leads to a new Weil-cohomology
whereas the second allows a comparison with classical $l$-adic cohomology.

\vspace{\abstand}

The plan of the paper is as follows:
In the second chapter, we begin by introducing ``$\hat{S}$-small'' categories for a superstructure $\hat{S}$
(categories whose morphisms form a set which is an element of $\hat{S}$)
and functors between these.
Given an enlargement $\hat{S}\rightarrow\widehat{\str{S}}$, we describe how to ``enlarge'' $\hat{S}$-small categories
and functors between them: We
assign an $\widehat{\str{S}}$-small category $\str{\cC}$ to every $\hat{S}$-small category $\cC$
and a functor $\str{F}:\str{\cC}\rightarrow\str{\cD}$ to every functor $F:\cC\rightarrow\cD$ of
$\hat{S}$-small categories.
Then we show that most basic properties of given $\hat{S}$-small categories and functors between them,
like the existence of certain finite limits or colimits, the representability of certain functors 
or the adjointness of a given pair of functors,
carry over to the associated enlargements.

In the third chapter, 
we study the enlargements of filtered and cofiltered categories
and the enlargements of filtered limits and colimits.
Using the ``saturation principle'' of enlargements,
we will be able to prove that all such limits in an $\hat{S}$-small category $\cC$
are ``dominated'' (in a sense that will be made precise) by objects in $\str{\cC}$.

In the fourth chapter, we will study $\hat{S}$-small additive and abelian
categories (and additive functors between them),
and their enlargements (which will be additive respectively abelian again), 
and we will give an interpretation of the enlargement of a category of $R$-modules 
(for a ring $R$ that is an element of $\hat{S}$)
as a category of \emph{internal} $\str{R}$-modules.

The fifth chapter is devoted to derived functors between $\hat{S}$-small abelian categories.
We will see that taking the enlargement is compatible with taking the $i$-th derived functor (again in a sense
that will be made precise).

In the sixth chapter, we look at triangulated and derived $\hat{S}$-small categories and exact functors between them
and study their enlargements. 
The enlargement of a triangulated category will be triangulated again,
the enlargement of a derived functor will be exact,
and we will prove several compatibilities between the various constructions.

In the seventh chapter, we add yet more structure and study $\hat{S}$-small fibred categories
and $\hat{S}$-small additive, abelian and triangulated fibrations and their enlargements.
This is important for most applications we have in mind.

Finally, we have added an appendix on basic definitions and properties of superstructures and enlargements for the convenience
of the reader.

\vspace{\abstand}

\def\Abe{\mathrm{Ab}}
\newcommand\MorC[1]{{\mathrm{Mor}}_{#1}}
\newcommand\Modfin[1]{{\mathrm{Mod}_{#1}^{\mathrm{fin}}}}
\def\image{\text{im}\,}
\def\kernel{\text{ker}\,}
\def\cokernel{\text{coker}\,}
\newcommand\rmod[1]{(#1-\underline{\mathrm{Mod}})}
\def\mR{\mathfrak{R}}
\def\Ensigns{\underline{\mathrm{Ens}}}

\section{Enlargements of categories}

\vspace{\abstand}

Let $\hat{S}$ be a superstructure and $*:\hat{S}\rightarrow\widehat{\str{S}}$
an enlargement. Even though a superstructure is not a universe in the sense of
\cite{SGA4I} (see \ref{bem:universe}), we use the same terminology as in
\cite{SGA4I} by simply replacing "universe" by "superstructure" in all
definitions. In particular, by an $\hat{S}$-small category $\cC$, we mean a
small category $\cC$ whose set of morphisms is contained in $\hat{S}$, and we
want to consider its image $\str{\cC}$ in $\widehat{\str{S}}$. In this
paragraph, we want to establish basic properties of $\str{\cC}$: It is again a
category (a $\widehat{\str{S}}$-small category to be precise), and it will
inherit a lot of the properties of $\cC$ (like having an initial object or
final object).

Not surprisingly, the proofs will rely heavily on the transfer principle, and therefore it will be convenient to
have short formal descriptions of the concepts involved. In particular, it will turn out to be easier to work
with the following definition of an $\hat{S}$-small category instead of the one given above:

\vspace{\abstand}

\begin{defi}\label{deficat}
  Let $\hat{S}$ be a superstructure.
  An \emph{$\hat{S}$-small category} is a quadruple
  $\langle M,s,t,c\rangle$ with
  \begin{enumerate}
    \item
      $M\in\hat{S}\setminus S$,
    \item\label{axiomst}
      $s,t:M\rightarrow M$ satisfying
      \begin{enumerate}
        \item
          $ss=ts=s$,
        \item
          $st=tt=t$,
      \end{enumerate}
    \item\label{axiomc}
      $c\subseteq M\times M\times M$ satisfying
      \begin{enumerate}
        \item
          $\forall f,g,h\in M:(\langle f,g,h\rangle\in c)\Rightarrow
          (sf=tg)\wedge(tf=th)\wedge(sg=sh)$,
        \item
          $\forall f,g\in M:(sf=tg)\Rightarrow(\exists!h\in M:\langle f,g,h\rangle\in c)$,
      \end{enumerate}
    \item\label{defcativ}
      $\forall f\in M:\langle f,tf,f\rangle,\langle sf,f,f\rangle\in c$,
    \item
      $
        \forall f_1,f_2,f_3,f_{12},f_{23},f_{123}\in M:
        \langle f_1,f_2,f_{12}\rangle,
        \langle f_{12},f_3,f_{123}\rangle,
        \langle f_2,f_3,f_{23}\rangle\in c
        \Rightarrow\langle f_1,f_{23},f_{123}\rangle\in c
      $.
  \end{enumerate}
  If $\cC=\langle M,s,t,c\rangle$ is an $\hat{S}$-small category, we call
  $\MorC{\cC}:=M$ the \emph{set of morphisms of $\cC$}.
\end{defi}

\vspace{\abstand}

\begin{bem}
Given an $\hat{S}$-small category $\cC=\langle M,s,t,c\rangle$ in the sense just defined,
we can get an $\hat{S}$-small category in the usual sense of \cite[I.1.0]{SGA4I} as follows:
\begin{enumerate}
  \item
    Define the set of \emph{objects} of $\cC$ as $\Ob(\cC):=s(M)$
    (note that also $\Ob(\cC)=t(M)$ because of \ref{deficat}\ref{axiomst};
    note also that because of \ref{deficat}\ref{defcativ}, elements of $s(M)=t(M)$ are units under composition,
    so objects correspond their identity morphisms).
  \item
    For objects $X,Y\in\Ob(\cC)$, define the set of \emph{morphisms from $X$ to $Y$}
    as $\Mor{\cC}{X}{Y}:=\{f\in M\vert sf=X\wedge tf=y\}$.
  \item
    For an object $X\in\Ob(\cC)$ define $\id{X}:=X\stackrel{\ref{deficat}\ref{axiomst}}{\in}\Mor{\cC}{X}{X}$.
  \item
    Consider objects $X,Y,Z\in\Ob(\cC)$ and morphisms $f\in\Mor{\cC}{Y}{Z}$ and $g\in\Mor{\cC}{X}{Y}$.
    According to \ref{deficat}\ref{axiomc}, there is a unique $h\in\Mor{\cC}{X}{Z}$ such that
    $\langle f,g,h\rangle\in c$. Put $fg:=f\circ g:=h$.
\end{enumerate}

\noindent
If on the other hand we are given an $\hat{S}$-small category $\cD$ in the sense of \cite{SGA4I} and if we put
\begin{itemize}
  \item
    $M:=\bigsqcup_{X,Y\in\Ob(\cD)}\Mor{\cD}{X}{Y}$,
  \item
    $sf:=\id{X}$ and $tf:=\id{Y}$ for $f\in\Mor{\cD}{X}{Y}$ and
  \item
    $c:=\bigl\{\langle f,g,h\rangle\,\vert\,\exists X,Y,Z\in\Ob(\cD):f\in\Mor{\cD}{Y}{Z}\wedge
    g\in\Mor{\cD}{X}{Y}\wedge f\circ g=h\bigr\}$,
\end{itemize}
then $\langle M,s,t,c\rangle$ is an $\hat{S}$-small category in the sense of \ref{deficat}.
\end{bem}

\vspace{\abstand}

\begin{bsp}\label{bspssmall}
  \mbox{}\\
  \begin{enumerate}
    \item\label{bspsmallcat}
      Let $\cC$ be a \emph{small} category. Then \ref{satzUsmallSsmall} shows that there is a superstructure
      $\hat{S}$ such that $\cC$ is isomorphic to an $\hat{S}$-small category.
    \item\label{bspensfin}
      Let $\hat{S}$ be a superstructure.
      There is an $\hat{S}$-small category $\cC$ that is equivalent to the category $\Ensfin$ of finite sets.
    \item\label{bsprmodfingen}
      Let $\hat{S}$ be a superstructure.
      If $A\in\hat{S}$ is a ring, then there is an $\hat{S}$-small category $\cC$ that is equivalent
      to $\Modfin{A}$, the category of finitely generated $A$-modules.
  \end{enumerate}
\end{bsp}

\vspace{\abstand}

As a next step, we want to describe the notion of a \emph{covariant functor} between
$\hat{S}$-small categories that are given in terms of \ref{deficat}:

\begin{defi}\label{defifunct}
  Let $\hat{S}$ be a superstructure,
  and let $\cC=\langle M,s,t,c\rangle$ and $\cD=\langle M',s',t',c'\rangle$ be $\hat{S}$-small categories.
  A \emph{(covariant) functor $F:\cC\rightarrow\cD$} is a map $F:M\rightarrow M'$ satisfying
  \begin{itemize}
    \item
      $Fs=s'F$,
    \item
      $Ft=t'F$ and
    \item
      $\forall f,g,h\in M:\langle f,g,h\rangle\in c\Rightarrow\langle Ff,Fg,Fh\rangle\in c'$.
  \end{itemize}
\end{defi}

\vspace{\abstand}

\begin{lemma}\label{lemmasmall}
  Let $\hat{S}$ be a superstructure, and let $\cC$ and $\cD$ be $\hat{S}$-small categories.
  Then the following categories are also $\hat{S}$-small:
  \begin{enumerate}
    \item\label{oppsmall}
      the opposite category $\cC^\circ$ of $\cC$,
    \item\label{productsmall}
      the product category $\cC\times\cD$,
    \item\label{covsmall}
      the category $\cov{\cC}{\cD}$ of covariant functors from $\cC$ to $\cD$,
    \item\label{consmall}
      the category $\con{\cC}{\cD}:=\cov{\cC^\circ}{\cD}$ of contravariant functors from $\cC$ to $\cD$.
  \end{enumerate}
\end{lemma}

\vspace{\abstand}

\begin{proof}
  Let $\cC=\langle M,s,t,c\rangle$ and $\cD=\langle M',s',t',c'\rangle$.
  Then \ref{oppsmall} and \ref{productsmall} are trivial, because $\cC^\circ=\langle M,t,s,c^\circ\rangle$
  with $\langle f,g,h\rangle\in c^\circ:\Leftrightarrow\langle g,f,h\rangle\in c$
  and $\cC\times\cD=\langle M\times M',s\times s',t\times t',c\times c'\rangle$.
  Because of
  \ref{defifunct}, the set of objects of $\cov{\cC}{\cD}$ is a subset of ${M'}^M$, and for two functors
  $F,G:\cC\rightarrow\cD$, a morphism $f$ between $F$ and $G$ is given by a family
  $\{f_X:FX\rightarrow GX\}_{X\in\Ob(\cC)}$ of morphisms in $\cD$, so that
  $\Mor{\cov{\cC}{\cD}}{F}{G}$ also is a subset of ${M'}^M$. Combining this, we see that the set of all morphisms
  in $\cov{\cC}{\cD}$ can be considered as a subset of
  $({M'}^M)^{{M'}^M\times{M'}^M}$ and is therefore an element of $\hat{S}$ which proves \ref{covsmall}.
  Finally, \ref{consmall} follows immediately from \ref{oppsmall} and \ref{covsmall}.
\end{proof}

\vspace{\abstand}

\noindent
Now we can look at the effect of applying $*$ to an $\hat{S}$-small category $\cC$:

\begin{satz}\label{satzstarcat}
  Let $*:\hat{S}\rightarrow\widehat{\str{S}}$ be an enlargement.
  and let $\cC$ be an $\hat{S}$-small category.
  Then $\str{\cC}$ is a $\widehat{\str{S}}$-small category with $\MorC{\str{\cC}}=\str{\MorC{\cC}}$
  (so in particular all morphisms in $\str{\cC}$ are internal).
\end{satz}

\vspace{\abstand}

\begin{proof}
  If $\cC=\langle M,s,t,c\rangle$,
  then $\str{\cC}=\langle\str{M},\str{s},\str{t},\str{c}\rangle$,
  and an easy transfer immediately proves that this new quadruple satisfies
  all the conditions in \ref{deficat}.
\end{proof}

\vspace{\abstand}

Before we go on, we need a technical lemma to compute the transfer of various formulas
from $\hat{S}$ to $\widehat{\str{S}}$
(by $\varphi[X_1,\ldots,X_n]$ we always mean a formula whose free variables are exactly $X_1,\ldots,X_n$):

\begin{lemma}\label{trst}
  In the situation of \ref{satzstarcat}, we have:
  \begin{enumerate}
    \item\label{trob}
      $\varphi[X]\equiv\bigl[X\in\Ob(\cC)\bigr]
      \;\;\Rightarrow\;\;
      \str{\varphi}[X]\equiv\bigl[X\in\Ob(\str{\cC})\bigr]$,
    \item\label{trmorph}
      $\varphi[f,X,Y]\equiv\bigl[f\in\Mor{\cC}{X}{Y}\bigr]
      \;\;\Rightarrow\;\;
      \str{\varphi}[f,X,Y]\equiv\bigl[f\in\Mor{\str{\cC}}{X}{Y}\bigr]$,
    \item\label{trcomp}
      $\varphi[f,g,h]\equiv\bigl[\text{$f$, $g$ composable morphisms in $\cC$ with $fg=h$}\bigr]$\\
      \mbox{}\hfill$\Rightarrow\;\;
      \str{\varphi}[f,g,h]\equiv\bigl[\text{$f$, $g$ composable morphisms in $\str{\cC}$ with $fg=h$}\bigr]$,
    \item\label{triso}
      $\varphi[f]\equiv\bigl[\text{$f$ isomorphism in $\cC$}\bigr]
      \;\;\Rightarrow\;\;
      \str{\varphi}[f]\equiv\bigl[\text{$f$ isomorphism in $\str{\cC}$}\bigr]$,
    \item\label{trmono}
      $\varphi[f]\equiv\bigl[\text{$f$ monomorphism in $\cC$}\bigr]
      \;\;\Rightarrow\;\;
      \str{\varphi}[f]\equiv\bigl[\text{$f$ monomorphism in $\str{\cC}$}\bigr]$,
    \item\label{trepi}
      $\varphi[f]\equiv\bigl[\text{$f$ epimorphism in $\cC$}\bigr]
      \;\;\Rightarrow\;\;
      \str{\varphi}[f]\equiv\bigl[\text{$f$ epimorphism in $\str{\cC}$}\bigr]$,
  \end{enumerate}
\end{lemma}

\vspace{\abstand}

\begin{proof}
  Let $\cC=\langle M,s,t,c\rangle$.\\[1 mm]
  \begin{enumerate}
    \item
      $\str{\varphi}[X]=\str{\bigl[\exists f\in M\,:\,sf=X\bigr]}
      =\bigl[\exists f\in\str{M}\,:\,(\str{s})f=X\bigr]
      =\bigl[X\in\Ob(\str{\cC})\bigr]$,\\[1 mm]
    \item
      $\str{\varphi}[f,X,Y]
       =\str{\bigl[(f\in M)\wedge(sf=X)\wedge(tf=Y)\bigr]}$\\
      \mbox{}\hfill
      $=\bigl[(f\in\str{M})\wedge((\str{s})f=X)\wedge(\str{t}f=Y)\bigr]
       =\bigl[f\in\Mor{\str{\cC}}{X}{Y}\bigr]$,
      \\[1 mm]
    \item
      $\str{\varphi}[f,g,h]=\str{\bigl[\langle f,g,h\rangle\in c\bigr]}=\bigl[\langle f,g,h\rangle\in\str{c}\bigr]$,
      \\[1 mm]
    \item
      $\str{\varphi}[f]=\str{\Bigl[\exists g\in M:\,\exists X,Y\in\Ob(\cC):\,(\langle f,g,X\rangle\in c)
      \wedge(\langle g,f,Y\rangle\in c
      )\Bigr]}$\\
      \mbox{}\hfill
      $\stackrel{\ref{trob}}{=}\Bigl[\exists g\in\str{M}:\,\exists X,Y\in\Ob(\str{\cC}):\,
      (\langle f,g,X\rangle\in\str{c})\wedge(\langle g,f,Y\rangle\in\str{c}))\Bigr]$,\\[1 mm]
    \item
      $
          \varphi[f]=\Bigl[f\in M\Bigr]\wedge
            \Bigl[\forall g_1,g_2\in M\,:\,
            \bigl(
              \exists h\in M\,:\,(\langle f,g_1,h\rangle\in c)\wedge(\langle f,g_2,h\rangle\in c)
            \bigr)
            \Rightarrow(g_1=g_2)
          \Bigr],
      $
      and from this the claim follows immediately.\\[1 mm]
    \item
      Analogous to ``monomorphism''!\\[1 mm]
  \end{enumerate}
\end{proof}

\vspace{\abstand}

\begin{cor}\label{corlemma}
  Let $*:\hat{S}\rightarrow\widehat{\str{S}}$ an enlargement,
  and let $\cC$ and $\cD$ be $\hat{S}$-small categories. Then we have:
  \begin{enumerate}
    \item\label{corobmor}
      $\Ob(\str{\cC})=\str{\Ob(\cC)}$, and
      if $m:\Ob(\cC)\times\Ob(\cC)\rightarrow\MorC{\cC}$ denotes the map
      that sends $\langle X,Y\rangle$ to $\Mor{\cC}{X}{Y}$,
      then $\str{m}:\Ob(\str{\cC})\times\Ob(\str{\cC})\rightarrow\MorC{\str{\cC}}$
      maps $\langle X,Y\rangle$ to $\Mor{\str{\cC}}{X}{Y}$.
      In particular we have
      $\Mor{\str{\cC}}{\str{X}}{\str{Y}}=\str{\Mor{\cC}{X}{Y}}$
      for $X,Y\in\Ob(\cC)$.
    \item\label{corstarfunct}
      $*$ induces a covariant functor from $\cC$ to $\str{\cC}$ that takes
      isomorphisms (resp. monomorphisms, resp. epimorphisms) to
      isomorphisms (resp. monomorphisms, resp. epimorphisms).
    \item\label{coropp}
      $\str{(\cC^\circ)}=(\str{\cC})^\circ$.
    \item\label{corproduct}
      $\str{(\cC\times\cD)}=\str{\cC}\times\str{\cD}$.
  \end{enumerate}
\end{cor}

\vspace{\abstand}

\begin{proof}
  Let $\cC=\langle M,s,t,c\rangle$ and $\cD=\langle M',s',t',c'\rangle$.
  \begin{enumerate}
    \item
      \ref{trst}\ref{trob},\ref{trmorph}.
    \item
      First of all, $*$ maps objects in $\cC$ to objects in $\str{\cC}$ and
      morphisms in $\cC$ to morphisms in $\str{\cC}$ because of
      \ref{trst}\ref{trob} resp. \ref{trmorph}, and it is compatible with compositions because of
      \ref{trst}\ref{trcomp}. Finally, if $X\in\Ob(\cC)$, then
      $\str{\id{X}}=\str{X}=\id{\str{X}}$.
    \item
      According to the proof of \ref{lemmasmall}\ref{oppsmall},
      $\cC^\circ=\langle M,t,s,c^\circ\rangle$ with
      $\langle f,g,h\rangle\in c^\circ\Leftrightarrow\langle g,f,h\rangle\in c$,
      and from this the claim follows immediately by transfer.
    \item
      According to the proof of \ref{lemmasmall}\ref{productsmall},
      $\cC^\circ=\langle M\times M',s\times s',t\times t',c\times c'\rangle$ from which the claim
      again immediately follows by transfer.
  \end{enumerate}
\end{proof}

\vspace{\abstand}

\begin{cor}\label{corisomorphic}
  Let $*:\hat{S}\rightarrow\widehat{\str{S}}$ be an enlargement,
  let $\cC$ be an $\hat{S}$-small category,
  and let $X$ and $Y$ be two objects of $\cC$.
  Then $X$ and $Y$ are isomorphic if and only if $\str{X}$ and $\str{Y}$ are isomorphic in $\str{\cC}$.
\end{cor}

\vspace{\abstand}

\begin{proof}
  We have
  \begin{multline*}
    \bigl[X\cong Y\bigr]
    \Leftrightarrow\bigl[\exists f\in\Mor{\cC}{X}{Y}:\text{$f$ isomorphism}\bigr] \\
    \stackrel{\text{\tiny transfer, \ref{trst}\ref{triso}, \ref{corlemma}\ref{corobmor}}}{\Longleftrightarrow}
    \bigl[\exists f\in\Mor{\str{\cC}}{\str{X}}{\str{Y}}:\text{$f$ isomorphism}\bigr]
    \Leftrightarrow\bigl[\str{X}\cong\str{Y}\bigr].
  \end{multline*}
\end{proof}

\vspace{\abstand}

\begin{defi}
  Let $\hat{S}$ be a superstructure. A \emph{finite} $\hat{S}$-small category is an $\hat{S}$-small category
  $\cC$ whose set  morphisms $\MorC{\cC}$ is a finite set.
\end{defi}

\vspace{\abstand}

For the next proposition,
recall from \ref{defifunct} that if $*:\hat{S}\rightarrow\widehat{\str{S}}$ is an enlargement
and if $\cC$ and $\cD$ are $\str{S}$-small categories,
then a functor from $\cC$ to $\cD$ is a map between the sets $\MorC{\cC}$ and $\MorC{\cD}$,
which are elements of $\widehat{\str{S}}$,
so in particular such a functor is an element of the superstructure $\widehat{\str{S}}$.
Therefore we can ask --- as for any element of $\widehat{\str{S}}$ --- whether a given functor is an \emph{internal}
element of $\widehat{\str{S}}$ (see \ref{defienlargement}).
Such functors, that are internal elements of $\widehat{\str{S}}$, we will call \emph{internal functors}.

\begin{satz}\label{satzfunctors}
  Let $*:\hat{S}\rightarrow\widehat{\str{S}}$ be an enlargement,
  and let $\cC$ and $\cD$ be $\hat{S}$-small categories.
  \begin{enumerate}
    \item\label{sfunctorsi}
      $\str{\cov{\cC}{\cD}}$ (resp. $\str{\con{\cC}{\cD}}$)
      is the subcategory of internal functors and internal morphisms of functors of
      $\cov{\str{\cC}}{\str{\cD}}$ (resp. $\con{\str{\cC}}{\str{\cD}}$).
      In particular, $*$ maps covariant (resp. contravariant) functors from $\cC$ to $\cD$
      to covariant (resp. contravariant) functors from $\str{\cC}$ to $\str{\cD}$.
    \item\label{sfunctorsiii}
      If $F:\cC\rightarrow\cD$ is a covariant functor, then the following diagram of categories and functors
      commutes:
      \[
        \xymatrix{
          {\cC} \ar[r]^{F} \ar[d]_{*} & {\cD} \ar[d]^{*} \\
          {\str{\cC}} \ar[r]_{\str{F}} & {\str{\cD}.} \\
        }
      \]
    \item\label{sfunctorsii}
      If $\cC$ is a \emph{finite} $\hat{S}$-small category, then
      \[
        \str{\cov{\cC}{\cD}}=\cov{\str{\cC}}{\str{\cD}}
        \;\;\;\;\;\;\text{and}\;\;\;\;\;\;
        \str{\con{\cC}{\cD}}=\con{\str{\cC}}{\str{\cD}}.
      \]
  \end{enumerate}
\end{satz}

\vspace{\abstand}

\begin{proof}
  Let $\cC=\langle M,s,t,c\rangle$ and $\cD=\langle M',s',t',c'\rangle$.
  We have the following statement in $\hat{S}$:
  \[
    \Ob(\cov{\cC}{\cD})=\{F:M\rightarrow M'\,\vert\,\text{$F$ subject to \ref{defifunct}
    concerning $\cC$, $\cD$}\}.
  \]
  Transfer of this gives us:
  \begin{multline*}
    \Ob(\str{\cov{\cC}{\cD}})
    \stackrel{\ref{corlemma}\ref{corobmor}}{=}
    \str{\Ob(\cov{\cC}{\cD})} \\
    =
    \{F:\str{M}\rightarrow\str{M'}\,\vert\,
    \text{$F$ \emph{internal} and subject to \ref{defifunct} concerning $\str{\cC}$, $\str{\cD}$}\},
  \end{multline*}
  i.e. the objects of $\str{\cov{\cC}{\cD}}$ are just the internal functors from $\str{\cC}$ to $\str{\cD}$.
  Now for the morphisms of functors, we have the following statement in $\hat{S}$:
  \begin{multline*}
    \forall F,G\in\Ob(\cov{\cC}{\cD})\,:\,
    \Mor{\cov{\cC}{\cD}}{F}{G} \\
    =\{
      f:\Ob(\cC)\rightarrow M'\,\vert\,\text{$f$ fulfills certain properties concerning $\cC$ and $\cD$}
    \},
  \end{multline*}
  and transferring this by using \ref{trst}\ref{trmorph} and \ref{corlemma}\ref{corobmor}, we get
  \begin{multline*}
    \forall F,G\in\Ob(\str{\cov{\cC}{\cD}})\,:\,
    \Mor{\str{\cov{\cC}{\cD}}}{F}{G} \\
    =\{
      f:\Ob(\str{\cC})\rightarrow\str{M'}\,\vert\,\text{$f$ \emph{internal},
      fulfills certain properties concerning $\str{\cC}$ and $\str{\cD}$}
    \},
  \end{multline*}
  i.e. for internal functors $F$ and $G$ from $\str{\cC}$ to $\str{\cD}$,
  the morphisms between $F$ and $G$ in $\str{\cov{\cC}{\cD}}$ are precisely the
  internal morphisms of functors between $F$ and $G$.
  This proves \ref{sfunctorsi} for covariant functors;
  the proof for contravariant functors is analogous.\\[2 mm]
  To prove \ref{sfunctorsiii}, we see immediately that for an arbitrary morphism $f\in\MorC{\cC}$ we have
  \[
    \str{F}\bigl[*(f)\bigr]=
    \str{F}\bigl[\str{f}\bigr]\stackrel{\ref{satzapp}\ref{satzappviii}}{=}
    \str{\bigl[F(f)\bigr]}=
    *\bigl[F(f)\bigr].
  \]
  To prove \ref{sfunctorsii}, according to \ref{sfunctorsi} we have to show that all functors from
  $\str{\cC}$ to $\str{\cD}$ and all morphisms of functors from $\str{\cC}$ to $\str{\cD}$
  are internal if $\cC$ is finite. But in this case, $*:\cC\rightarrow\str{\cC}$ is a \emph{bijection},
  i.e. both functors from $\str{\cC}$ to $\str{\cD}$ and morphisms of functors from $\str{\cC}$ to $\str{\cD}$
  are maps from a \emph{finite} set to an internal set.
  But such maps are internal (see \ref{satzinternal}\ref{satzinternalfinite}).
\end{proof}

\vspace{\abstand}

\begin{cor}
  Let $*:\hat{S}\rightarrow\widehat{\str{S}}$ an enlargement,
  and let $\cC$ be an $\hat{S}$-small category.
  \begin{enumerate}
    \item\label{corsubcat}
      Let $\cD$ be a (full) subcategory of $\cC$. Then $\cD$ is $\hat{S}$-small,
      and $\str{\cD}$ is a $\widehat{\str{S}}$-small (full) subcategory of $\str{\cC}$.
    \item\label{corCX}
      Let $X$ be an object of $\cC$, then $\cC/X$, the category of objects over $X$,
      is $\hat{S}$-small,
      and $\str{(\cC/X)}=\str{\cC}/\str{X}$.
    \item\label{corfunres}
      Let $\cD$ be a subcategory of $\cC$, and let $F$ be a covariant (resp. contravariant)
      functor from $\cC$ into another $\hat{S}$-small category $\cC'$. Then
      $\str{(F\vert\cD)}=\str{F}\vert\str{\cD}:\str{\cD}\rightarrow\str{\cC'}$.
    \item\label{corsieve}
      Let $X$ be an object of $\cC$,
      and let $\cR$ be a sieve of $X$ (see \cite[I.4.1]{SGA4I}),
      i.e. a full subcategory of $\cC/X$ with the property that if $Y$ is an object of $\cC/X$
      and if there exists a morphism from $Y$ to an object of $\cR$ in $\cC/X$,
      then $Y$ belongs to $\cR$.
      Then $\str{\cR}$ is a sieve of $\str{X}$.
  \end{enumerate}
\end{cor}

\vspace{\abstand}

\begin{proof}
  Let $\cC=\langle M,s,t,c\rangle$,
  and let $\cD$ be a subcategory of $\cC$. Then by definition we have
  $\Ob(\cD)\subseteq\Ob(\cC)$,
  $\forall X,Y\in\Ob(\cD):\Mor{\cD}{X}{Y}\subseteq\Mor{\cC}{X}{Y}$,
  and that the inclusion $\cD\hookrightarrow\cC$ is a \emph{functor}.
  Obviously, these conditions are equivalent to
  $\cD=\langle N,u,v,d\rangle$ with $N\subseteq M$, $u=s\cap(N\times N)$,
  $v=t\cap(N\times N)$ and $d=c\cap(N\times N\times N)$. This immediately shows that $\cD$ is
  $\hat{S}$-small, and transfer proves that $\str{\cD}$ is a subcategory of $\str{\cC}$.\\[1mm]
  If $\cD$ is a \emph{full} subcategory of $\cC$,
  then we have the true statement $\forall X,Y\in\Ob(\cD):\Mor{\cD}{X}{Y}=\Mor{\cC}{X}{Y}$ in $\hat{S}$,
  and transfer of this (using \ref{corlemma}\ref{corobmor}) shows that $\str{\cD}$ is a full subcategory of $\str{\cC}$ ---
  so the proof of \ref{corsubcat} is complete.\\[2mm]
  Let $X$ be an object of $\cC$.
  Put
  \[
    \begin{array}{l}
      \displaystyle N:=\Bigl\{\langle f,g,h\rangle\in c\Bigl\vert t(g)=t(h)=X\Bigr\}, \\[3mm]
      \displaystyle u:N\rightarrow N,\ \langle f,g,h\rangle\mapsto s(f), \\[2mm]
      \displaystyle v:N\rightarrow N,\ \langle f,g,h\rangle\mapsto t(f), \\[2mm]
      \displaystyle d:=\Bigl\{\langle f_1,g_1,h_1\rangle,\langle f_2,g_2,h_2\rangle,\langle f,g,h\rangle
        \in N\times N\times N\Bigl\vert
        \langle f_1,f_2,f\rangle\in c
      \Bigr\}.
    \end{array}
  \]
  Then $\cC/X=\langle N,u,v,d\rangle$, i.e. $\cC/X$ is $\hat{S}$-small.
  The following statement is true in $\hat{S}$:
  \begin{multline*}
    \Bigl[\Ob(\cC/X)=\bigl\{\langle Y,f\rangle\in\Ob(\cC)\times M\bigl\vert f\in\Mor{\cC}{Y}{X}\bigr\}\Bigr] \\
    \wedge\Bigl[\forall\langle Y,f\rangle,\langle Z,g\rangle\in\Ob(\cC/X):
      \Mor{\cC/X}{\langle Y,f\rangle}{\langle Z,g\rangle}
      =\bigl\{h\in\Mor{\cC}{Y}{Z}\bigl\vert g\circ h=f\bigr\}
    \Bigr].
  \end{multline*}
  Using \ref{corlemma}\ref{corobmor}, we see that transfer of this proves \ref{corCX}.\\[2mm]
  \ref{corfunres} is trivial, following immediately by transfer from the statement
  $\forall X\in\Ob(\cD):(F\vert\cD)(X)=F(X)$
  (resp. the corresponding statement concerning morphisms).\\[2 mm]
  Finally, let $\cR$ be a sieve of an object $X$ of $\cC$.
  Then \ref{corsubcat} and \ref{corCX} show that $\str{\cR}$ is a full subcategory of $\str{\cC}/\str{X}$.
  In addition to that, we have the following true statement in $\hat{S}$:
  \[
    \forall Y\in\Ob(\cC/X):
    \forall R\in\Ob(\cR):
    \bigl[\Mor{\cC/X}{Y}{R}\neq\emptyset\bigr]
    \Rightarrow
    \bigl[Y\in\Ob(\cR)\bigr].
  \]
  By transfer, we get the following true statement in $\widehat{\str{S}}$:
  \[
    \forall Y\in\underbrace{\str{\Ob(\cC/X)}}_{
      \stackrel{\ref{corCX}\ref{corlemma}\ref{corobmor}}{=}\Ob(\str{\cC}/\str{X})}:
    \forall R\in\underbrace{\str{\Ob(\cR)}}_{
      \stackrel{\ref{corlemma}\ref{corobmor}}{=}\Ob(\str{\cR})}:
    \bigl[\underbrace{\str{\Mor{\cC/X}{Y}{R}}}_{
      \stackrel{\ref{corCX}\ref{corlemma}\ref{corobmor}}{=}\Mor{\str{\cC}/\str{X}}{Y}{R}
    }\neq\emptyset\bigr]
    \Rightarrow
    \bigl[Y\in\underbrace{\str{\Ob(\cR)}}_{
      \stackrel{\ref{corlemma}\ref{corobmor}}{=}\Ob(\str{\cR})
    }\bigr],
  \]
  i.e. $\str{\cR}$ indeed is a sieve of $\str{X}$.
\end{proof}

\vspace{\abstand}

\begin{defi}
  For a set of sets $N$, let $\Ens{N}$ denote the category of sets that are elements of $N$.
\end{defi}

\vspace{\abstand}

\begin{bem}
  Note that if $N\in\hat{S}$ is a set of sets
  in a superstructure $\hat{S}$,
  then $\Ens{N}$ is an $\hat{S}$-small
  category (because the set of morphisms in $\Ens{N}$ is then obviously an element of $\hat{S}$).
  Furthermore, if $\hat{S}$ is an arbitrary superstructure,
  then any full subcategory of the category of sets that is $\hat{S}$-small equals
  $\Ens{N}$ for a suitable $N\in\hat{S}$.
\end{bem}

\vspace{\abstand}

\begin{satz}\label{satzens}
  Let $*:\hat{S}\rightarrow\widehat{\str{S}}$ be an enlargement,
  and let $N\in\hat{S}$ be a set of sets.
  Then $\str{\Ens{N}}$ is the subcategory of $\Ens{\str{N}}$ that consists of all the objects but
  whose morphisms are only the \emph{internal} maps.
\end{satz}

\vspace{\abstand}

\begin{proof}
  Transfer of the statement
  $\Ob\Ens{N}\stackrel{\ref{deficat}}{=}\left\{\id{X}\vert X\in N\right\}\cong N$
  gives us $\Ob(\str{\Ens{N}})\cong\str{N}$,
  so the objects of $\str{\Ens{N}}$ are exactly the objects of $\Ens{\str{N}}$. As for the morphisms,
  transfer of
  \[\forall \id{X},\id{Y}\in\Ob\Ens{N}\,:\,\Mor{\Ens{N}}{\id{X}}{\id{Y}}=Y^X\]
  gives the desired result.
\end{proof}

\vspace{\abstand}

\begin{cor}\label{correp}
  Let $*:\hat{S}\rightarrow\widehat{\str{S}}$ be an enlargement,
  let $N\in\hat{S}$ be a set of sets,
  and let $\cC$ be an $\hat{S}$-small category with $M:=\MorC{\cC}\subseteq N$.
  \begin{enumerate}
    \item\label{correpi}
      Let $F:\cC^\circ\rightarrow\Ens{N}$ be a functor
      (i.e. $F$ is a presheaf on $\cC$ with values in $\Ens{N}$).
      Then $\str{F}$ is a presheaf on $\str{\cC}$ with values in $\Ens{\str{N}}$.
    \item\label{correph}
      If in particular we look at the presheaf
      $h_X:=\Mor{\cC}{\_}{X}:\cC^\circ\rightarrow\Ens{M}\subseteq\Ens{N}$ for an element $X\in\Ob(\cC)$,
      then
      $\str{h_X}=\Mor{\str{\cC}}{\_}{\str{X}}=:h_{\str{X}}$.
    \item\label{correpii}
      Consider the statement
      \[
        \varphi[F,X]\equiv\Bigl[\bigl(F\in\Ob(\con{\cC}{\Ens{N}})\bigr)\wedge\bigl(X\in\Ob(\cC)\bigr)\wedge
        \bigl(\text{$F$ is representable by $X$}\bigr)\Bigl]
      \]
      (note that this makes sense because of $M\subseteq N$). Then \\
      \begin{minipage}{15cm}\begin{multline*}
        \str{\varphi}[F,X]=
          \Bigl[
            \text{$F$ is an internal presheaf on $\str{\cC}$ with values in $\Ens{\str{N}}$} \\
            \text{which is internally representable by $X$}
          \Bigr].
      \end{multline*}\end{minipage}\\[1 mm]
      (By ``internally representable'' we mean that
      there is an object $Y\in\Ob(\str{\cC})$ and an \emph{internal} isomorphism of functors
      between $F$ and
      $\Mor{\str{\cC}}{\_}{Y}$.)
    \item\label{correpiii}
      If a presheaf $F$ on $\cC$ with values in $\Ens{N}$ is representable by an object $X\in\Ob(\cC)$,
      then $\str{F}$ is representable by $\str{X}$.
    \item\label{correpiv}
      If all presheaves on $\cC$ with values in $\Ens{\str{N}}$ are representable, then all \emph{internal}
      presheaves on $\str{\cC}$ with values in $\Ens{\str{N}}$ are representable.
  \end{enumerate}
\end{cor}

\vspace{\abstand}

\begin{proof}
  \ref{satzfunctors}\ref{sfunctorsi} and \ref{satzens} show that $\str{F}$ is a functor from
  $\str{\cC}$ to $\Ens{\str{M}}$ which is a subcategory of $\Ens{\str{N}}$,
  so this proves \ref{correpi}.\\[1 mm]
  For \ref{correph}, we just have to transfer the statement
  \[
    \forall Y\in\Ob(\cC):
    h_X(Y)=\Mor{\cC}{Y}{X}
  \]
  by using \ref{corlemma}\ref{corobmor}.\\[1 mm]
  To prove \ref{correpii},
  let $h:\cC\rightarrow\con{\cC}{\Ens{N}},X\mapsto h_X$ denote the canonical fully faithful embedding
  where
  --- as in \ref{correph} --- $h_X$ denotes the functor $\Mor{\cC}{\_}{X}:\cC^\circ\rightarrow\Ens{N}$.
  Then
  \begin{multline*}
    \varphi[F,X]=\Bigl[
      \bigl(F\in\Ob(\con{\cC}{\Ens{N}})\bigr)\wedge\bigl(X\in\Ob(\cC)\bigr) \\
      \wedge\bigl(\exists f\in\Mor{\con{\cC}{\Ens{N}}}{F}{h(X)}
      \,:\,
      \text{$f$ isomorphism in $\con{\cC}{\Ens{N}}$}\bigr)
    \Bigr].
  \end{multline*}
  We then get
  \begin{multline*}
    \str{\varphi}[F,X]\stackrel{\ref{trst}\ref{triso},\ref{corlemma}\ref{corobmor}}{=}\Bigl[
      \bigl(F\in\Ob(\str{\con{\cC}{\Ens{N}}})\bigr)\wedge\bigl(X\in\Ob(\str{\cC})\bigr) \\
      \wedge\bigl(\exists f\in\Mor{\str{\con{\cC}{\Ens{N}}}}{F}{(\str{h})(X)}\,:\,
      \text{$f$ isomorphism in $\str{\con{\cC}{\Ens{N}}}$}\bigr)
    \Bigr]
  \end{multline*}
  and by \ref{satzfunctors}\ref{sfunctorsi}
  \begin{multline*}
    \str{\varphi}[F,X]=\Bigl[
      \bigl(
        \text{$F$ is an internal presheaf on $\cC$ with values in $\Ens{\str{N}}$}
      \bigr)
      \wedge\bigl(
        X\in\Ob(\str{C})
      \bigr) \\
      \wedge\bigl(
        \text{$F$ and $(\str{h})(X)$ are internally isomorphic}
      \bigr)
    \Bigr].
  \end{multline*}
  What is $(\str{h})(X)$?
  \begin{multline*}
    \forall X,Y\in\Ob(\cC)\,:\,[h(X)](Y)=\Mor{\cC}{Y}{X} \\
    \stackrel{\text{\tiny transfer}, \ref{trst}\ref{trmorph}, \ref{corlemma}\ref{corobmor}}
      {\Longleftrightarrow}
    \forall X,Y\in\Ob(\str{\cC})\,:\,[(\str{h})(X)](Y)=\Mor{\str{\cC}}{Y}{X}=h_{X}(Y). \\
  \end{multline*}
  So we get:
  \begin{multline*}
    \str{\varphi}[F,X]=
      \Bigl[
        \text{$F$ is an internal presheaf on $\cC$ with values in $\Ens{\str{N}}$} \\
        \text{which is internally representable by $X$}
      \Bigr],
  \end{multline*}
  which is \ref{correpii}.\\[2 mm]
  If a given $F_0\in\Ob(\con{\cC}{\Ens{N}})$ is representable by $X_0\in\Ob(\cC)$, then
  $\varphi[F_0,X_0]$ is true, and by transfer $\str{\varphi}[\str{F_0},\str{X_0}]$ must be true,
  i.e. $\str{F_0}$ is (internally) representable by $\str{X_0}$, which proves \ref{correpiii}.\\[2 mm]
  If finally all such $F$ are representable, then the statement
  \[
    \forall F\in\Ob(\con{\cC}{\Ens{N}})\,:\,
    \exists X\in\Ob(\cC)\,:\,
    \varphi[F,X]
  \]
  is true, and by transfer we get that for all internal presheaves $F$ on $\str{\cC}$ with values
  in $\Ens{\str{N}}$, there exists a $X\in\Ob(\str{\cC})$ that (internally) represents them.
\end{proof}

\vspace{\abstand}

\begin{cor}\label{corfinlim}
  Let $*:\hat{S}\rightarrow\widehat{\str{S}}$ be an enlargement,
  let $\cI$ be a \emph{finite} $\hat{S}$-small category,
  and let $\cC$ be an $\hat{S}$-small category.
  If for all functors $G:\cI\rightarrow\cC$ the projective limit $\varprojlim G$
  (resp. inductive limit $\varinjlim G$) exists in $\cC$, then
  for all functors $G:\cI\rightarrow\str{\cC}$, $\varprojlim G$
  (resp. $\varinjlim G$) exists in $\str{\cC}$.
\end{cor}

\vspace{\abstand}

\begin{proof}
  Let $k:\Ob(\cC)\rightarrow\Ob(\cov{\cI}{\cC})$, $X\mapsto k_X$ denote the map that maps $X$ to
  the constant functor with value $X$.
  By assumption, we have the following statement in $\hat{S}$:
  \[
    \forall G\in\Ob(\cov{\cI}{\cC})\,:\,
    \text{$\left\{X\mapsto\Mor{\cov{\cI}{\cC}}{k_X}{G}\right\}_{X\in\Ob(\cC)}$ is representable}.
  \]
  According to \ref{lemmasmall}\ref{covsmall}, the category $\cov{\cI}{\cC}$ is $\hat{S}$-small,
  i.e. of the form $\langle M,s,t,c\rangle$ for a $M\in\hat{S}$. This implies that the presheaf
  $X\mapsto\Mor{\cov{\cI}{\cC}}{c_X}{G}$ on $\cC$ has values in $\Ens{M}$ for any functor
  $G:\cI\rightarrow\cC$,
  so we can apply \ref{correp}\ref{correpii} to get by transfer:
  \[
    \forall G\in\Ob(\underbrace{\str{\cov{\cI}{\cC}}}_
      {\stackrel{\ref{satzfunctors}\ref{sfunctorsii}}{=}\cov{\str{\cI}}{\str{\cC}}}
    )\,:\,
    \text{$\Bigl\{X\mapsto\Mor{\underbrace{\str{\cov{\cI}{\cC}}}_
      {\stackrel{\ref{satzfunctors}\ref{sfunctorsii}}{=}\cov{\str{\cI}}{\str{\cC}}}
    }{(\str{k})(X)}{G}\Bigr\}_{X\in\Ob(\str{\cC})}$ is representable}.
  \]
  It is easy to see that $(\str{k})(X)$ is again the constant functor with value $X$, so we have proved
  that
  $\varprojlim G$ exists in $\str{\cC}$
  for all functors $G$ from $\str{\cI}$ to $\str{\cC}$.
  But as we noted in the proof of \ref{satzfunctors}\ref{sfunctorsii}, the map $\cI\xrightarrow{*}\str{\cI}$
  is a \emph{bijection} (and thus because of \ref{corlemma}\ref{corstarfunct} an isomorphism of categories),
  so that we have a canonical bijection between functors from $\str{\cI}$ to $\str{\cC}$ and
  functors from $\cI$ to $\str{\cC}$, which proves the corollary for projective limits.
  The case of inductive limits is analogous.
\end{proof}

\vspace{\abstand}

\begin{cor}\label{corfinlims}
  Let $*:\hat{S}\rightarrow\widehat{\str{S}}$ be an enlargement,
  and let $\cC$ be an $\hat{S}$-small category.
  Then if $\cC$ has one of the following properties, then so has $\str{\cC}$:
  \begin{enumerate}
    \item
      $\cC$ has an initial object.
    \item
      $\cC$ has a final object.
    \item
      $\cC$ has a null object.
    \item
      Arbitrary finite direct sums exist in $\cC$.
    \item
      Arbitrary finite direct products exist in $\cC$.
    \item
      Arbitrary finite fibred sums exist in $\cC$.
    \item
      Arbitrary finite fibred products exist in $\cC$.
    \item
      Difference cokernels of two arbitrary morphisms exist in $\cC$.
    \item
      Difference kernels of two arbitrary morphisms exist in $\cC$.
  \end{enumerate}
\end{cor}

\vspace{\abstand}

\begin{proof}
  These are all special cases of \ref{corfinlim}, because all the properties can be expressed in terms
  of finite projective or inductive limits.
\end{proof}

\vspace{\abstand}

\begin{satz}\label{satzadjoint}
  Let $*:\hat{S}\rightarrow\widehat{\str{S}}$ be an enlargement,
  let $\cC$ and $\cD$ be $\hat{S}$-small categories,
  and let $F:\cC\rightarrow\cD$ and $G:\cD\rightarrow\cC$ be two covariant functors.
  If $F$ is left adjoint to $G$, then $\str{F}:\str{\cC}\rightarrow\str{\cD}$
  is left adjoint to $\str{G}:\str{\cD}\rightarrow\str{\cC}$.
\end{satz}

\vspace{\abstand}

\begin{proof}
  First of all, according to \ref{satzfunctors}\ref{sfunctorsi}, $\str{F}$ is a functor
  from $\str{\cC}$ to $\str{\cD}$,
  and $\str{G}$ is a functor from $\str{\cD}$ to $\str{\cC}$.\\[1mm]
  Now choose a set $N\in\hat{S}$ that contains the (disjoint) union of the set of morphisms in $\cC$
  and the set of morphisms in $\cD$, and consider the following two covariant functors
  $\alpha,\beta:\cC^\circ\times\cD\rightarrow\Ens{N}$:
  \[
    \xymatrix@R=2mm@C=5cm{
      {\cC^\circ\times\cD} \ar@<1mm>[r] \ar@<-1mm>[r] & {\Ens{N}} \\
      & {\Mor{\cD}{FX}{Y}} \\
      {(X,Y)} \ar@{|->}[ur]^{\alpha} \ar@{|->}[dr]_{\beta} & \\
      & {\Mor{\cC}{X}{GY}.} \\
    }
  \]
  Saying that $F$ is left adjoint to $G$ is equivalent to saying that $\alpha$ and $\beta$ are
  isomorphic in the category $\cov{\cC^\circ\times\cD}{\Ens{N}}$ which is
  $\hat{S}$-small because of \ref{lemmasmall}\ref{oppsmall}, \ref{productsmall} and \ref{covsmall}.
  Because of $\ref{corisomorphic}$, this is equivalent to $\str{\alpha}$ an $\str{\beta}$ being isomorphic
  in $\str{\cov{\cC^\circ\times\cD}{\Ens{N}}}$.
  Now we can apply
  \ref{corlemma}\ref{coropp} and \ref{corproduct},
  \ref{satzfunctors}\ref{sfunctorsi},
  and \ref{satzens}
  to see that this is a subcategory of
  $\cov{(\str{\cC})^\circ\times\str{\cD}}{\Ens{\str{N}}}$,
  i.e. if $F$ is left adjoint to $G$, then $\str{\alpha}$ and $\str{\beta}$ are isomorphic in this latter
  category, and it is clear by transfer and by \ref{corlemma}\ref{corobmor} that $\str{\alpha}$ and
  $\str{\beta}$ are given as follows:
  \[
    \xymatrix@R=2mm@C=5cm{
      {(\str{\cC})^\circ\times\str{\cD}} \ar@<1mm>[r] \ar@<-1mm>[r] & {\Ens{\str{N}}} \\
      & {\Mor{\str{\cD}}{(\str{F})X}{Y}} \\
      {(X,Y)} \ar@{|->}[ur]^{\str{\alpha}} \ar@{|->}[dr]_{\str{\beta}} & \\
      & {\Mor{\str{\cC}}{X}{(\str{G})Y},} \\
    }
  \]
  i.e. the fact that $\str{\alpha}$ and $\str{\beta}$ are isomorphic is equivalent to $\str{F}$ being
  left adjoint to $\str{G}$.
\end{proof}

\vspace{\abstand}

\section{Limits and enlargements}

\vspace{\abstand}

So far, so good! --- But until now, the only property of enlargements we have used is the transfer principle,
and therefore we have not revealed any properties of enlarged categories $\str{\cC}$ that the original category
$\cC$ has not already had. This will change in the next proposition.\\[2 mm]

\noindent
Remember that a category $\cI$ is called \emph{pseudo-cofiltered} if
the following two conditions are satisfied (compare \cite[I.2.7]{SGA4I}):
\begin{itemize}
  \item
    For any pair of morphisms $\varphi_1:i_1\rightarrow i$ and
    $\varphi_2:i_2\rightarrow i$,
    there are morphisms $\psi_1:j\rightarrow i_1$ and
    $\psi_2:j\rightarrow i_2$ in $\cI$ such that
    $\varphi_1\circ\psi_1=\varphi_2\circ\psi_2$:
    \xymatrix@R=0mm{
      & {i_1} \ar[dr]^{\varphi_1} & \\
      {j} \ar@{.>}[ur]^{\exists\psi_1} \ar@{.>}[dr]_{\exists\psi_2} & & {i} \\
      & {i_2} \ar[ur]_{\varphi_2} & \\
    }
  \item
    For any pair of morphisms $\varphi_1,\varphi_2:i_1\rightarrow i_2$ in $\cI$,
    there is a morphism $\psi:j\rightarrow i_1$ in $\cI$ such that
    $\varphi_1\circ\psi=\varphi_2\circ\psi$:
    \xymatrix{
      {j} \ar@{.>}[r]^{\exists\psi} & {i_1} \ar@<1mm>[r]^{\varphi_1} \ar@<-1mm>[r]_{\varphi_2} & {i_2} \\
    }
\end{itemize}

\noindent
$\cI$ is called \emph{cofiltered} if it is pseudo-cofiltered, not empty and connected,
where in this case being connected is equivalent to the condition that for any pair of objects
$i_1,i_2\in\Ob(\cI)$, there is a pair of morphisms $\psi_1:j\rightarrow i_1$ and $\psi_2:j\rightarrow i_2$:
\xymatrix@R=0mm{
  & {i_1} \\
  {j} \ar@{.>}[ur]^{\exists\psi_1} \ar@{.>}[dr]_{\exists\psi_2} & \\
  & {i_2} \\
}\\[2 mm]
A category $\cI$ is called \emph{pseudo-filtered}, if $\cI^\circ$ is pseudo-cofiltered,
and $\cI$ is called \emph{filtered} if $\cI^\circ$ is cofiltered.

\vspace{\abstand}

\begin{satz}\label{satziinfty}
  Let $*:\hat{S}\rightarrow\widehat{\str{S}}$ be an enlargement,
  and let $\cI$ be an $\hat{S}$-small category.
  \begin{enumerate}
    \item\label{iinftycofil}
      If $\cI$ is cofiltered, then there exists an object $i_{-\infty}\in\Ob(\str{\cI})$,
      together with
      a family of morphisms ${\{p_i:i_{-\infty}\rightarrow\str{i}\}}_{i\in\Ob(\cI)}$ with the property that
      for all morphisms $i_1\xrightarrow{\varphi}i_2$ in $\cI$,
      the following triangle of morphisms in $\str{\cI}$
      commutes:
      \xymatrix@R=0mm{
        & {\str{i_1}} \ar[dd]^{\str{\varphi}} \\
        {i_{-\infty}} \ar[ru]^{p_{i_1}} \ar[rd]_{p_{i_2}} & \\
        & {\str{i_2}} \\
      }
    \item\label{iinftyfil}
      If $\cI$ is filtered, then there exists an object $i_\infty\in\Ob(\str{\cI})$,
      together with
      a family of morphisms ${\{\iota_i:\str{i}\rightarrow i_\infty\}}_{i\in\Ob(\cI)}$ with the property that
      for all morphisms $i_1\xrightarrow{\varphi}i_2$ in $\cI$,
      the following triangle of morphisms in $\str{\cI}$
      commutes:
      \xymatrix@R=0mm{
        {i_1} \ar[dd]_{\str{\varphi}} \ar[rd]^{\iota_{i_1}} & \\
        & {i_\infty} \\
        {i_2} \ar[ru]_{\iota_{i_2}} & \\
      }
  \end{enumerate}
\end{satz}

\vspace{\abstand}

\begin{proof}
  We only give a proof for \ref{iinftycofil}, the proof for \ref{iinftyfil} follows immediately by looking at
  the opposite category $\cI^\circ$.\\[1 mm]
  We start by proving the following statement:
  \begin{equation}\label{eqfinsys}
    \text{\begin{minipage}{14cm}
    Let $\cJ$ be a finite subsystem of $\cI$. By this we mean a selection of finitely many objects of $\cI$
    and of finitely many morphisms between those objects
    (note that $\cJ$ in general will not be a subcategory
    of $\cI$). Then there exists an object $i_\cJ$ of $\cI$ and a family of morphisms
    $p_j^{\cJ}:i_\cJ\rightarrow j$ for every object $j$ of $\cI$ that is in $\cJ$ such that for
    any morphism $j_1\xrightarrow{\varphi}j_2$ contained in $\cJ$,
    the following triangle of morphisms in $\cI$ commutes:
    \xymatrix@R=0mm{
      & {j_1} \ar[dd]^{\varphi} \\
      {i_\cJ} \ar[ru]^{p_{j_1}^{\cJ}} \ar[rd]_{p_{j_2}^{\cJ}} & \\
      & {j_2} \\
    }\end{minipage}
    }
  \end{equation}
  We prove \eqref{eqfinsys} by induction over the number of objects and morphisms in $\cJ$:\\[2 mm]
  If $\cJ$ is empty, then we can take any object of $\cI$ as $i_\cJ$
  (note that by definition a cofiltered category
  is not empty).\\[2 mm]
  If $\cJ$ contains exactly one object $j$ and no morphisms, then we can take $i_\cJ:=j$ and
  $p_j^{\cJ}:=\id{j}$.\\[2 mm]
  Let us assume that we have proven \eqref{eqfinsys} for systems of at most $n$ objects
  that contain no morphisms at all,
  and let $\cJ$ be a system with $(n+1)$ objects with no morphisms at all. Let $j_0$ be one of these objects,
  and let $\cJ':=\{j_1,\ldots,j_n\}$ be the system of the other objects.
  By induction, we have an object $i_{\cJ'}$ of $\cI$ and morphisms $p_{j_k}^{\cJ'}:i_{\cJ'}\rightarrow j_k$
  for $k=1,\ldots,n$.
  Because $\cI$ is connected,
  we find an object $i_{\cJ}\in\Ob(\cI)$ and morphisms
  $i_{\cJ}\xrightarrow{p_{j_0}^{\cJ}}j_0$ and $i_{\cJ}\xrightarrow{q}i_{\cJ'}$:
  \xymatrix@R=0mm@C=15mm{
    & {j_0} \\
    {i_\cJ} \ar[ur]^{p_{j_0}^{\cJ}} \ar[dr]_{q} & \\
    & {i_{\cJ'}} \\
  }\\
  Put $p_{j_k}^{\cJ}:=p_{j_k}^{\cJ'}\circ q$ for $k=1,\ldots,n$:
  \xymatrix@R=0mm@C=20mm{
    & & {j_0} \\
    {i_\cJ} \ar[urr]^{p_{j_0}^{\cJ}} \ar[ddr]_{q} & & \\
    & & {j_1} \\
    & {i_{\cJ'}} \ar[ur]^{p_{j_1}^{\cJ'}} \ar[dr]_{p_{j_n}^{\cJ'}} & {\vdots} \\
    & & {j_n} \\
  }\\
  Finally, let us assume we have proven \eqref{eqfinsys} for systems with arbitrarily
  (but finitely) many objects that
  contain at most $n$ morphisms, and let $\cJ$ be a system with $n+1$ morphisms.
  Let $j_0\xrightarrow{\varphi_0}j_0'$ be one of the morphisms,
  and let $\cJ'$ be the system that contains all the objects of $\cJ$ and
  all the morphisms except $\varphi_0$. Then by induction, we have $i_{\cJ'}$ and morphisms
  $p_j^{\cJ'}:i_{\cJ'}\rightarrow j$ for every object $j$ contained in $\cJ$, and all triangles
  \xymatrix@R=0mm{
    & {j_1} \ar[dd]^{\varphi} \\
    {i_{\cJ'}} \ar[ru]^{p_{j_1}^{\cJ}} \ar[rd]_{p_{j_2}^{\cJ}} & \\
    & {j_2} \\
  }
  commute for morphisms $j_1\xrightarrow{\varphi}j_2$ in $\cJ'$, but possibly not for $\varphi_0$.
  Because $\cI$ is pseudo-cofiltered, we find an object $i_\cJ\in\Ob(\cI)$
  and a morphism $i_\cJ\xrightarrow{q}i_{\cJ'}$ with
  $p_{j_0'}^{\cJ'}\circ q=(\varphi_0\circ p_{j_0'}^{\cJ'})\circ q$:
  \xymatrix@R=0mm{
    {i_\cJ} \ar[r]^{q} & {i_{\cJ'}} \ar@(ur,ul)[rr]^{p_{j_0'}^{\cJ'}}
      \ar[dr]_{p_{j_0}^{\cJ'}} & & {j_0'} \\
    & & {j_0} \ar[ur]_{\varphi_0} & \\
  }\\
  Now we just set $p_j^{\cJ}:=p_j^{\cJ}\circ q$ for all objects $j$ in $\cJ$
  and thereby prove \eqref{eqfinsys} in general.\\[2 mm]
  Let $M\in\hat{S}$ denote the set of all morphisms in $\cI$,
  and let $F$ denote the set of all finite subsystems of $\cI$ --- we certainly have $F\in\hat{S}$.
  For $\cJ\in F$, define a set $U_\cJ$ as follows:
  \begin{multline*}
    U_\cJ=\Bigl\{
      \langle i,p\rangle\,\Bigl\vert
      i\in\Ob(\cI),
      p\in M^{\Ob(\cI)}, \\
      \forall j\in\{\text{objects in $\cJ$}\}:p(j)\in\Mor{\cI}{i}{j},
      \forall j_1\xrightarrow{\varphi}j_2\in\{\text{morphisms in $\cJ$}\}:p(j_2)=\varphi\circ
      p(j_1)
    \Bigr\}
  \end{multline*}
  Statement \eqref{eqfinsys} shows two things:
  \begin{itemize}
    \item
      For all $\cJ\in F$, the set $U_\cJ$ is not empty (note that for objects $i$ not contained in $\cJ$,
      $p(i)$ can be chosen arbitrarily in $M$, we can for example set $p(i):=\id{i}$).
    \item
      For $\cJ_1,\ldots,\cJ_n\in F$, the intersection $U_{\cJ_1}\cap\ldots\cap U_{\cJ_n}$ is not empty
      (apply \eqref{eqfinsys} to $\cJ:=\cJ_1\cup\ldots\cup\cJ_n$).
  \end{itemize}
  According to the saturation principle \ref{defienlargement}\ref{saturation},
  these imply that the intersection
  $\bigcap_{\cJ\in F}\str{U_\cJ}$ is not empty,
  i.e. we find a pair $\langle i_{-\infty},p\rangle$, consisting of an object $i_{-\infty}\in\Ob(\str{\cI})$
  and an internal map $p:\Ob(\str{\cI})\rightarrow\str{M}$ in this intersection.
  For an object $i$ of $\cI$, there certainly is a $\cJ\in F$ containing $i$, and because
  $\langle i_{-\infty},p\rangle$ is in $\str{U_\cJ}$, we conclude that $p(\str{i})$ is a morphism
  from $i_{-\infty}$ to $\str{i}$ and can thus set $p_i:=p(\str{i})$.\\[1 mm]
  If $i_1\xrightarrow{\varphi}i_2$ is any morphism in $\cI$, we find a $\cJ\in F$ that contains
  $i_1$, $i_2$ and $\varphi$,
  and the fact that $\langle i_{-\infty},p\rangle$ is an element of $\str{U_\cJ}$ implies that the
  associated triangle commutes. This concludes the proof of the proposition.
\end{proof}

\vspace{\abstand}

\begin{bsp}(Construction of algebraic closures)\\
  Let $k$ be a field. To construct an algebraic closure of $k$, it suffices to construct an extension $K$ of $k$
  such that $K$ contains all finite algebraic extensions of $k$, because then taking the algebraic closure of $k$
  in $K$ yields an algebraic closure of $k$.

  To construct such a field $K$,
  choose a set of sets $U$ such that for every finite algebraic extension $L$ of $k$,
  there is a bijection from the set which underlies $L$ to an element of $U$.

  Consider the category $\cC$ whose objects are pairs $\langle L,s\rangle$ with $L$ an element of $U$
  and $s$ a field-structure on $L$ turning $L$ into a finite algebraic extension of $k$,
  and whose morphisms are defined as
  \[
    \Mor{\cC}{\langle L_1,s_1\rangle}{\langle L_2,s_2\rangle}:=
    \left\{\begin{array}{cl}
      \{*\} & \text{if there is a $k$-embedding $L_1\hookrightarrow L_2$,} \\
      \emptyset & \text{otherwise.}
    \end{array}\right.
  \]
  We then can find a base set $S$ such that $\cC$ is (isomorphic to) an $\hat{S}$-small category
  and such that $U$ is an element of $\hat{S}\setminus S$ (compare \ref{bspssmall}\ref{bspsmallcat}).

  It is easy to see that $\cC$ is \emph{filtered}, so that
  we find an object $i_\infty$ of $\str{\cC}$ and morphisms
  $\iota_L:\str{L}\rightarrow i_\infty$ for every object $L$ of $\cC$ as in \ref{satziinfty}\ref{iinftyfil}.

  By transfer, $i_\infty$ defines a pair $\langle K,s\rangle$ where $K$ is a set (and an element of $\str{U}$)
  and $s$ is a field-structure on $K$.
  The existence of the $\iota_L$ and transfer show
  that we have $k$-embeddings $L\hookrightarrow K$ for all finite algebraic extensions $L$ of $k$,
  and we are done.\\[1 mm]

  Note that in order to prove the existence of enlargements (compare \ref{thmenlargement}),
  methods similar to the construction of algebraic closures in elementary algebra are used;
  therefore it is not really surprising that we are able to give a short proof for the existence of algebraic closures
  as soon as we have the powerful tool of enlargements at our disposal.
\end{bsp}

\vspace{\abstand}

\begin{cor}\label{corphi}
  Let $*:\hat{S}\rightarrow\widehat{\str{S}}$ be an enlargement,
  let $\cI$ be an $\hat{S}$-small category,
  and let $\varphi[X]$ be a formula in $\hat{S}$.
  \begin{enumerate}
    \item\label{corcofilphi}
      Assume that $\cI$ is cofiltered,
      let $i_{-\infty}$ and $\{p_i\}$ be as in \ref{satziinfty}\ref{iinftycofil},
      and assume that for all morphisms $\psi:i\rightarrow j$ in $\cI$,
      statement $\varphi[i]$ implies $\varphi[j]$. Then
      \[
        \Bigl(\forall i\in\Ob(\cI):\varphi[i]\Bigr)
        \;\;\;\Longleftrightarrow\;\;\;
        \str{\varphi}[i_{-\infty}].
      \]
    \item\label{corfilphi}
      Assume that $\cI$ is filtered,
      let $i_\infty$ and $\{\iota_i\}$ be as in \ref{satziinfty}\ref{iinftyfil},
      and assume that for all morphisms $\psi:i\rightarrow j$ in $\cI$,
      statement $\varphi[j]$ implies $\varphi[i]$. Then
      \[
        \Bigl(\forall i\in\Ob(\cI):\varphi[i]\Bigr)
        \;\;\;\Longleftrightarrow\;\;\;
        \str{\varphi}[i_\infty].
      \]
  \end{enumerate}
\end{cor}

\vspace{\abstand}

\begin{proof}
  We only have to prove \ref{corcofilphi}, because \ref{corfilphi} follows from this by looking at
  $\cI^\circ$.
  If $\varphi[i]$ is true for all $i\in\Ob(\cI)$, then by transfer and by \ref{corlemma}\ref{corobmor},
  $\str{\varphi}[i]$ is true for
  all $i\in\Ob(\str{\cI})$, so it is in particular true for $i:=i_{-\infty}\in\Ob(\str{\cI})$.
  On the other hand, by assumption, the following statement is true in $\hat{S}$:
  \[
    \forall j,k\in\Ob(\cI):
    \Bigl[\Mor{\cI}{j}{k}\neq\emptyset\Bigr]
    \Rightarrow
    \Bigl[\varphi[j]\Rightarrow\varphi[k]\Bigr],
  \]
  and transfer of this gives us
  \[
    \forall j,k\in\Ob(\str{\cI}):
    \Bigl[\Mor{\str{\cI}}{j}{k}\neq\emptyset\Bigr]
    \Rightarrow
    \Bigl[\str{\varphi}[j]\Rightarrow\str{\varphi}[k]\Bigr].
  \]
  If we specialize this to $j:=i_{-\infty}$ and $k:=\str{i}$ for an $i\in\Ob(\cI)$ and note that
  $\Mor{\str{\cI}}{i_{-\infty}}{\str{i}}$ is not empty (because it contains $p_i$),
  we get $\bigl[\str{\varphi}[i_{-\infty}]\Rightarrow\str{\varphi}[\str{i}]\bigr]$.
  But by transfer, the statements $\str{\varphi}[\str{i}]$ and $\varphi[i]$ are equivalent,
  and this completes the proof.
\end{proof}

\vspace{\abstand}

\begin{bsp}\label{bsp:garbenbed}
  Let $X$ be a topological space,
  let $U\subseteq X$ be an open subspace,
  and let $\cF$ be a presheaf on $X$.

  Recall that a \emph{covering} of $U$ is a family $\cU=(U_j)_{j\in\cJ}$ of open subspaces of $U$
  with $\bigcup_{j\in\cJ}U_j=U$,
  and a \emph{refinement} of $\cU$ is a covering $\cV=(V_k)_{k\in\cK}$ of $U$,
  such that there exists a map
  $f:\cK\rightarrow\cJ$ satisfying $V_k\subseteq U_{f(k)}$ for all $k\in\cK$.
  Consider the following partially ordered set $\cI$:
  An element of $\cI$ is represented by a covering of $U$,
  and two coverings define the same element if each is a refinement of the other.
  The partial order is defined by
  \[
    \cV\leq\cU:\Leftrightarrow\text{$\cV$ is a refinement of $\cU$.}
  \]
  If $\cU$ and $\cV$ are coverings of $U$, then
  $\cU\cap\cV:=(U_j\cap V_k)_{(j,k)\in\cJ\times\cK}$
  is obviously a covering of $U$ and a refinement of both $\cU$ and $\cV$,
  which shows that $\cI$, considered as a category, is cofiltered.

  For a covering $\cU$ of $U$,
  recall the \emph{sheaf condition for $\cF$ with respect to $\cU$},
  which states that the following sequence is exact:
  \[
    \cF(U)\rightarrow\prod_{j\in\cJ}\cF(U_j)\rightrightarrows\prod_{i,j\in\cJ}\cF(U_i\cap U_j).
  \]
  It is easy to see that if $\cV$ is a refinement of $\cU$
  and if $\cF$ satisfies the sheaf condition with respect to $\cV$,
  then it also satisfies the sheaf condition with respect to $\cU$.
  This in particular implies that if $\cU$ and $\cV$ are coverings of $U$ that define the same element of $\cI$,
  then $\cF$ satisfies the sheaf condition with respect to $\cU$ iff it satisfies the sheaf condition with respect
  to $\cV$. It therefore makes sense to say that $\cF$ satisfies the sheaf condition with respect to an element of $\cI$.

  Now let $\hat{S}$ be a superstructure that contains the small category $\cI$,
  and let $\varphi[X]$ be the following formula in $\hat{S}$:
  \[
    \varphi[X]\equiv\Bigl[\text{$X$ is an element of $\cI$,
      and $\cF$ satisfies the sheaf condition with respect to $X$}\Bigr].
  \]
  By proposition \ref{satziinfty}\ref{iinftycofil}, we find an object $i_{-\infty}$ of $\str{\cI}$,
  represented by a ``hyper covering" $\cU=(U_j)_{j\in\cJ}$ of $\str{U}$,
  where $\cU$ is an internal family of $\str{\text{open}}$ subsets of $\str{U}$ which covers $\str{U}$
  and refines all standard coverings of $\cU$.

  Applying corollary \ref{corphi}\ref{corcofilphi} to this situation, we get
  that $\cF$ satisfies the sheaf condition with respect to \emph{every covering} iff
  $\str\cF$ satisfies the $\str{\text{sheaf condition}}$ for the \emph{one hyper covering} $\cU$.
  This latter condition is obviously satisfied iff the following two conditions are both satisfied:
  \begin{enumerate}
    \item
      If $s,t\in\str{(\cF(U))}=(\str{\cF})(U)$ satisfy $s\vert U_j=t\vert U_j$ for all $j\in\cJ$, then $s=t$.
    \item
      If $(s_j)_{j\in\cJ}$ is an internal family with $s_j\in(\str{\cF})(U_j)$,
      satisfying $s_j\vert U_j\cap U_{j'}=s_{j'}\vert U_j\cap U_{j'}$ for all $j,j'\in\cJ$,
      then there exists an $s\in\str{(\cF(U))}$ with $s_j=s\vert U_j$ for all $j\in\cJ$.
  \end{enumerate}
\end{bsp}

\vspace{\abstand}

\begin{cor}\label{cormap}
  Let $*:\hat{S}\rightarrow\widehat{\str{S}}$ be an enlargement,
  let $\cI$ be an $\hat{S}$-small category,
  let $N\in\hat{S}$ be a set of sets,
  and let $F:\cI\rightarrow\Ens{N}$ be a covariant functor.
  \begin{enumerate}
    \item\label{cormapcofil}
      Assume that $\cI$ is cofiltered, and let $i_{-\infty}$ be as in
      \ref{satziinfty}\ref{iinftycofil}.
      Then we have a canonical injection $\mu:\varprojlim F\hookrightarrow\str{F}(i_{-\infty})$,
      and for every object $i\in\Ob(\cI)$, we have a commutative diagram
      \begin{equation}\label{eqmu}
        \xymatrix{
          & {\str{F}(i_{-\infty})} \ar[dr]^{\str{F}(p_i)} \\
          {\varprojlim F} \ar[rr]_{(x_i)_{i\in\Ob(\cI)}\mapsto\str{x_i}} \ar@{^{(}-->}[ur]^{\mu} & &
            {\str{F}(\str{i})} \\
        }
      \end{equation}
    \item\label{cormapfil}
      Assume that $\cI$ is filtered, and let $i_\infty$ and $\{\iota_i\}$ be as in
      \ref{satziinfty}\ref{iinftyfil}.
      Then we have a canonical injection $\nu:\varinjlim F\hookrightarrow\str{F}(i_\infty)$.
  \end{enumerate}
\end{cor}

\vspace{\abstand}

\begin{proof}
  In the situation of \ref{cormapcofil},
  let $x$ be an element of $\varprojlim F$, given by a family $(x_i)_{i\in\Ob(\cI)}$ of elements
  $x_i\in F(i)$ such that for all morphisms $\varphi:i\rightarrow j$ in $\cI$, the element $x_i$
  is mapped to $x_j$ under $F(\varphi)$.
  If we consider each $x_i$ as a morphism $\{*\}\rightarrow F(i)$ in $\Ens{N}$,
  we consequently can consider $x$ as a map $x:\Ob(\cI)\rightarrow\MorC{\Ens{N}}$ with the property
  \[
    \Bigl[\forall i\in\Ob(\cI):x(i)\in\Mor{\Ens{N}}{\{*\}}{Fi}\Bigr]
    \wedge
    \Bigl[\forall i,j\in\Ob(\cI):\forall\varphi\in\Mor{\cI}{i}{j}:x(j)=F\varphi\circ x(i)\Bigr].
  \]
  Then by transfer,
  $\str{x}$ is a map from
  $\Ob(\str{\cI})$ to
  $\MorC{\str{\Ens{N}}}\stackrel{\ref{satzens}}{\hookrightarrow}{\MorC{\Ens{\str{N}}}}$
  satisfying
  \begin{multline*}
    \Bigl[\forall i\in\Ob(\str{\cI}):x(i)\in\Mor{\Ens{\str{N}}}{\{*\}}{\str{F}(i)}\Bigr] \\
    \wedge
    \Bigl[\forall i,j\in\Ob(\str{\cI}):\forall\varphi\in\Mor{\str{\cI}}{i}{j}:x(j)=\str{F}(\varphi)\circ x(i)\Bigr],
  \end{multline*}
  i.e. by setting $\tilde{x}_i:=\bigl[\str{x}(i)\bigr](*)\in\str{F}(i)$ for all objects $i\in\Ob(\str{I})$,
  we get a compatible system $(\tilde{x}_i)_{i\in\Ob(\str{I})}$,
  and we define $\mu(x)$ as $\tilde{x}_{i_{-\infty}}$.\\[1 mm]
  To see that the so defined map $\mu$ is injective, let $y=(y_i)$ be another element of $P$ with the property
  $\mu(x)=\mu(y)$. We have to show that $x_i=y_i$ for all objects $i$ of $\cI$.
  Look at the formula $\varphi[X]\equiv\bigl[x_X=y_X\bigr]$ in $\hat{S}$.
  We have to show that the statement $\varphi[i]$ is true for all $i\in\Ob(\cI)$.
  If $\varphi:i\rightarrow j$ is a morphism in $\cI$,
  then $\varphi[i]$ implies $\varphi[j]$, so we can apply \ref{corphi}\ref{corcofilphi}, i.e.
  we have to show that $\str{\varphi}[i_{-\infty}]$ is true in $\widehat{\str{S}}$.
  But obviously $\str{\varphi}[X]=\bigl[\tilde{x}_X=\tilde{y}_X\bigr]$, so $\str{\varphi}[i_{-\infty}]$ is
  true because of $\mu(x)=\mu(y)$.
  So we see that $\mu$ indeed is injective.\\[2 mm]
  In the situation of \ref{cormapfil},
  let $x$ be an element of $\varinjlim F$, represented by an element $x_j$ of $F(j)$ for an $j\in\Ob(\cI)$.
  Put $\nu(x):=\iota_j(\str{x_j})$. To see that this is well defined, it suffices to show that for
  a morphism $\varphi:j\rightarrow k$ in $\cI$ we have $\iota_j(\str{x_j})=\iota_k\bigl[\str{(\varphi(x_j))}\bigr]$,
  but this is clear by definition of the $\{\iota_i\}$.
  So we get a well defined map $\nu:\varinjlim F\rightarrow F(i_\infty)$.\\[1 mm]
  To show that $\nu$ is injective, let $y$ be another element of $\varinjlim F$,
  represented by $y_k\in F(k)$, with $x\neq y$. We have to show $\nu(x)\neq\nu(y)$.
  Look at the formula
  \[
    \varphi[X]\equiv\Bigl[
      \forall\varphi\in\Mor{\cI}{j}{X}:
      \forall\psi\in\Mor{\cI}{k}{X}:
      \varphi(x_j)\neq\psi(y_k)
    \Bigr].
  \]
  Because we assume $x\neq y$, we know that $\varphi[i]$ is true for all $i\in\Ob(\cI)$.
  Then by transfer $\str{\varphi}[i]$ must be true for all $i\in\Ob(\str{\cI})$, so that in particular
  $\str{\varphi}[i_\infty]$ is true in $\widehat{\str{S}}$
  (to see this, we could have also applied \ref{corphi}\ref{corfilphi}). Now
  \[
    \varphi[i_\infty]=\Bigl[
      \forall\varphi\in\Mor{\cI}{\str{j}}{i_\infty}:
      \forall\psi\in\Mor{\cI}{\str{k}}{i_\infty}:
      \varphi(\str{x_j})\neq\psi(\str{y_k})
    \Bigr],
  \]
  and if we specialize this to $\varphi:=\iota_j$ and $\psi:=\iota_k$, we see that in particular
  $\nu(x)\neq\nu(y)$ is true, i.e. $\nu$ is indeed injective.
\end{proof}

\vspace{\abstand}

Let $\cI$ be a cofiltered category,
let $\cC$ be any category,
and let $G:\cI\rightarrow\cC$ be a covariant functor.
Remember that by definition, the projective limit $\varprojlim G$ is the presheaf
$X\mapsto\varprojlim_{i\in\Ob(\cI)}\Mor{\cC}{X}{Gi}$ on $\cC$,
and that we say that ``the projective limit exists in $\cC$"
if $\varprojlim G$ is representable by an object $P$ of $\cC$.

Assume for a moment that $\cI$ has an initial object $i_0$. In that case, $\varprojlim G$ obviously exists in $\cC$,
and it is represented by $Gi_0$, i.e. we have a commutative diagram of presheaves
\[
  \xymatrix{
    {\varprojlim G} \ar@{^{(}->}[d]_{\prod({\text{proj}_i})} \ar[r]^{\sim} &
      {h_{Gi_0}} \ar@{^{(}->}[d]^{\prod G({i_0\rightarrow i})} \\
    {\prod_{i\in\Ob(\cI)}h_{Gi}} \ar@{=}[r] & {\prod_{i\in\Ob(\cI)}h_{Gi}.}
  }
\]
Now, in general $\cI$ will not have an initial object, but we have seen in \ref{satziinfty}\ref{iinftycofil}
that in the category
$\str{\cI}$, there is an object $i_{-\infty}$ which is ``nearly as good as an initial object for $\cI$",
and we could hope that in the following commuting diagram of presheaves
\[
  \xymatrix{
    {\varprojlim G} \ar@{^{(}->}[d]_{\prod({\text{proj}_i})} \ar[r] &
      {\tilde{G}:=\bigl[X\mapsto\Mor{\str{\cC}}{X}{\str{G}i_{-\infty}}\bigr]} \ar@{^{(}->}[d]^{\prod{p_i}} \\
    {\prod_{i\in\Ob(\cI)}h_{Gi}} \ar[r]_-{*} &
    {\prod_{i\in\Ob(\cI)}\bigl[X\mapsto\Mor{\str{\cC}}{\str{X}}{\str{G}i}\bigr]}
  }
\]
the morphism $\varprojlim G\rightarrow\tilde{G}$ is an isomorphism.
Unfortunately, there are two problems with this:
First, for an object $X$ of $\cC$, an object $i$ of $\cI$ and a section $\varphi\in\tilde{G}(X)$,
${p_i}_*(\varphi)$ must be of the form $\str{\varphi_i}$ for a suitable $\varphi_i\in\Mor{\cC}{X}{Gi}$,
and there is no reason why this should be the case.

Second, remember that the pair $\langle i_{-\infty},(p_i)_i\rangle$ is in general not uniquely determined.
For example, for any morphism $p:i'_{-\infty}\rightarrow i_{-\infty}$ in $\str{\cI}$,
the pair $\langle i'_{-\infty},(p_i\circ p)_i\rangle$
has the same property.

Now let $\tilde{G}'$ be the presheaf $\bigl[X\mapsto\Mor{\str{\cC}}{X}{\str{G}i'_{-\infty}}\bigr]$,
and for an object $X$ of $\cC$, let $s$ and $t$ be two distinct sections in $\tilde{G}'(X)$.
If $p$ is not a monomorphism,
it can happen that $s$ and $t$ become equal in $\tilde{G}(X)$, and this shows that $\prod{p_i}_*$
in general will not be a monomorphism.

To solve the first problem, we will define a subpresheaf $\finG$ of $\tilde{G}$
whose sections are mapped to something of the right form under $\prod{p_i}_*$.
To solve the second problem, we will introduce an equivalence relation $\sim$ on $\finG$ such that
two section of $\tilde{G}$ that are mapped to the same section under $\prod{p_i}_*$
will be equivalent.

Having done that, we will be able to prove $\varprojlim G\cong\finG/\sim$.

\vspace{\abstand}

\begin{cor}\label{corinfinite}
  Let $*:\hat{S}\rightarrow\widehat{\str{S}}$ be an enlargement,
  and let $G:\cI\rightarrow\cC$ be a covariant functor of $\hat{S}$-small categories.
  \begin{enumerate}
    \item\label{corcofil}
      Assume that $\cI$ is cofiltered, and let $i_{-\infty}$ and $\{p_i\}$ be as in
      \ref{satziinfty}\ref{iinftycofil}.
      Define the presheaf $\tilde{G}$ on $\cC$ by
      $\tilde{G}(X):=\Mor{\str{\cC}}{\str{X}}{(\str{G})(i_{-\infty})}$ for $X\in\Ob(\cC)$,
      define the subpresheaf $\finG$ of $\tilde{G}$ by
      \[
        \finG(X):=\left\{\left.\varphi\in\tilde{G}(X)\right\vert
          \forall i\in\Ob(\cI):\exists\varphi_i\in\Mor{\cC}{X}{G(i)}:(\str{G})(p_i)\circ\varphi=\str{\varphi_i}
        \right\}
      \]
      for $X\in\Ob(\cC)$,
      and define an equivalence relation $\sim$ on $\finG$ as follows:
      \[
        \forall X\in\Ob(\cC):
        \forall\varphi,\psi\in\finG(X):
        \varphi\sim\psi:\Leftrightarrow
        \forall i\in\Ob(\cI):
        (\str{G})(p_i)\circ\varphi=(\str{G})(p_i)\circ\psi.
      \]
      Then we have a canonical monomorphism
      $\varprojlim G\hookrightarrow\tilde{G}$ of presheaves on $\cC$ that
      induces an isomorphism $\varprojlim G\xrightarrow{\sim}\finG/\sim$.
      If in addition $\Mor{\cC}{X}{G(i)}$ is \emph{finite} for all $X\in\Ob(\cC)$ and all $i\in\Ob(\cI)$,
      then $\finG=\tilde{G}$.\\[2 mm]
      The following diagram of presheaves on $\cC$ commutes:
      \begin{equation}\label{eqpresheavescolim}
        \xymatrix{
          {\varprojlim G} \ar@{^{(}->}[d]_{\prod({\text{proj}_i})} \ar[r]^-{\sim} &
            {\finG/\sim} \ar@{^{(}->}[d]^{\prod{p_i}} \\
          {\prod_{i\in\Ob(\cI)}h_{Gi}} \ar[r]_-{*} &
          {\prod_{i\in\Ob(\cI)}\bigl[X\mapsto\Mor{\str{\cC}}{\str{X}}{\str{G}i}\bigr].}
        }
      \end{equation}
    \item\label{corfil}
      Assume that $\cI$ is filtered, and let $i_\infty$ and $\{\iota_i\}$ be as in
      \ref{satziinfty}\ref{iinftyfil}.
      Define the presheaf $\tilde{G}$ on $\cC^\circ$ by
      $\tilde{G}(X):=\Mor{\str{\cC}}{(\str{G})(i_\infty)}{\str{X}}$ for $X\in\Ob(\cC)$,
      define the subpresheaf $\finG$ of $\tilde{G}$ by
      \[
        \finG(X):=\left\{\left.\varphi\in\tilde{G}(X)\right\vert
          \forall i\in\Ob(\cI):\exists\varphi_i\in\Mor{\cC}{G(i)}{X}:\varphi\circ(\str{G})(\iota_i)=\str{\varphi_i}
        \right\}
      \]
      for $X\in\Ob(\cC)$,
      and define an equivalence relation $\sim$ on $\finG$ as follows:
      \[
        \forall X\in\Ob(\cC):
        \forall\varphi,\psi\in\finG(X):
        \varphi\sim\psi:\Leftrightarrow
        \forall i\in\Ob(\cI):
        \varphi\circ(\str{G})(\iota_i)=\psi\circ(\str{G})(\iota_i).
      \]
      Then we have a canonical monomorphism
      $\varinjlim G\hookrightarrow\tilde{G}$ of presheaves on $\cC^\circ$ that
      induces an isomorphism $\varinjlim G\xrightarrow{\sim}\finG/\sim$.
      If in addition $\Mor{\cC}{G(i)}{X}$ is \emph{finite} for all $X\in\Ob(\cC)$ and all $i\in\Ob(\cI)$,
      then $\finG=\tilde{G}$.\\[2 mm]
      For every $i\in\Ob(\cI)$, the following diagram of presheaves on $\cC^\circ$ commutes:
      \begin{equation}\label{eqpreshaveslim}
        \xymatrix{
          {\varinjlim G} \ar@{^{(}->}[d]_{\prod({\text{proj}_i})} \ar[r]^-{\sim} &
            {\finG/\sim} \ar@{^{(}->}[d]^{\prod{\iota_i}} \\
          {\prod_{i\in\Ob(\cI)}h_{Gi}} \ar[r]_-{*} &
          {\prod_{i\in\Ob(\cI)}\bigl[X\mapsto\Mor{\str{\cC}}{\str{G}i}{\str{X}}\bigr].}
        }
      \end{equation}
  \end{enumerate}
\end{cor}

\vspace{\abstand}

\begin{proof}
  We only prove \ref{corcofil}, because
  \ref{corfil} follows from \ref{corcofil} by considering $\cI^\circ$ and $\cC^\circ$.
  Let $X$ be an object of $\cC$.
  The first claim is equivalent to the statement that for every object $X\in\Ob(\cC)$ there is a canonical injection
  \[
    \alpha_X:\varprojlim\Mor{\cC}{X}{G(i)}\hookrightarrow\Mor{\str{\cC}}{\str{X}}{(\str{G})(i_{-\infty})}
  \]
  which is functorial in $X$.
  To see this, apply \ref{cormap}\ref{corcofil} to $N:=\MorC{\cC}$ and
  $F:\cI\rightarrow\Ens{N}$, $i\mapsto\Mor{\cC}{X}{Gi}$
  (the functorality is clear by construction). We thus get a monomorphism $\varprojlim G\hookrightarrow\tilde{G}$
  of presheaves on $\cC$ which we call $\alpha$.\\[2 mm]
  We have already proved that the image of this monomorphism is contained in $\finG$ --- this is just \eqref{eqmu},
  so that we get an induced map $\varprojlim G\xrightarrow{\bar{\alpha}}\finG/\sim$,
  and \eqref{eqmu} obviously also shows that $\bar{\alpha}$ is injective.
  Furthermore, it is clear by construction that
  \eqref{eqpresheavescolim} indeed is a well-defined, commuting diagram of presheaves on $\cC$.\\[1 mm]
  It remains to be seen that $\bar{\alpha}$ is surjective.
  So let $p:\str{X}\rightarrow(\str{G})(i_{-\infty})$ be an element of $\finG(X)$.
  For $i\in\Ob(\cI)$, by definition of $\finG$, there exists a $\varphi_i:X\rightarrow G(i)$
  such that $(\str{G})(p_i)\circ p=\str{\varphi_i}$.
  If $i\xrightarrow{\psi}j$ is a morphism in $\cI$, then by definition of the $p_i$ we have
  \[
    \begin{array}{ll}
      & \psi\circ p_i=p_j \\[1 mm]
      \Rightarrow & (\str{G})(\str{\psi})\circ(\str{G})(p_i)=(\str{G})(p_j) \\[1 mm]
      \Rightarrow & (\str{G})(\str{\psi})\circ\bigl[(\str{G})(p_i)\circ p\bigr]=\bigl[(\str{G})(p_j)\circ p\bigr] \\[1 mm]
      \Rightarrow & (\str{G})(\str{\psi})(\str{\varphi_i})=\str{\varphi_j} \\[1 mm]
      \Leftrightarrow & G(\psi)(\varphi_i)=\varphi_j \\
    \end{array}
  \]
  which proves that $\varphi:=\{\varphi_i\}_{i\in\Ob(\cI)}$ defines an element of
  $(\varprojlim G)(X)$. By construction we have $\alpha_X(\varphi)\sim p$,
  i.e. $\varphi$ is a preimage of $p$ under $\bar{\alpha}$.
  This shows that $\bar{\alpha}$ indeed is an isomorphism of presheaves on $\cC$.\\[2 mm]
  Finally, let us assume that $\Mor{\cC}{X}{G(i)}$ is finite for all $X\in\Ob(\cC)$ and
  $i\in\Ob(\cI)$. Then $*$ induces bijections
  $\Mor{\cC}{X}{G(i)}\xrightarrow{\sim}\Mor{\str{\cC}}{\str{X}}{(\str{G})(\str{i})}$
  because of \ref{satzapp}\ref{satzappfin},
  i.e. all elements on the right hand side are of the form $\str{\varphi}$ for a suitable $\varphi$,
  so that all sections of $\tilde{G}$ belong to $\finG$.
\end{proof}

\vspace{\abstand}

\begin{bsp}\label{bsplim}
  \mbox{}\\
  \begin{enumerate}
    \item\label{bspprojlim}
      Let $\cC$ be the category of finite rings
      (whose elements are contained in the set $\Z$, say),
      let $\cI$ be the category of the  ordered set $\N_+$ (ordered by $\geq$),
      let $\hat{S}$ be a superstructure that contains $\cI$ and $\cC$ as $\hat{S}$-small categories,
      let $p$ be a prime number,
      let $G:\cI\rightarrow\cC$ be the functor $n\mapsto\Z/p^n\Z$,
      and let $*:\hat{S}\rightarrow\widehat{\str{S}}$ be an enlargement.\\[2 mm]
      Then $\cI$ is cofiltered, and we can apply \ref{corinfinite}\ref{corcofil} to $\varprojlim G$:
      $\cI$ is the category of the ordered set $\str{\N}$,
      as $i_{-\infty}$, we can choose any $h\in\str{\N}\setminus\N$, the category $\str{\cC}$ is the category of
      \strfin\ rings (with internal ring-structure, whose elements are contained in $\str{\Z}$)
      and internal ring-homomorphisms as morphisms,
      and $\str{G}$ is the functor that sends $h$ to $\str{\Z}/p^h\str{\Z}$.
      The corollary then tells us (because all the sets of morphisms in $\cC$ are finite) that to give
      a compatible system of ring-homomorphisms from a finite Ring $R$ into every $\Z/p^n\Z$ is the same
      as giving a ring-homomorphism (which is automatically internal because $R$ is finite) from $R$ into
      the \strfin\ ring $\str{\Z}/p^h\str{\Z}$ where two such ring-homomorphisms give the same
      compatible system if and only if their compositions with all projections
      $\str{\Z}/p^h\str{\Z}\twoheadrightarrow\Z/p^n\Z$ coincide. This is equivalent to the existence of
      a ring-isomorphism $\Z_p\xrightarrow{\sim}\bigl(\str{\Z}/p^h\str{\Z}\bigr)/I$
      where $I$ is the (external) ideal generated
      by $\{p^k\vert k\in\str{\N}\setminus\N\}$.\\[1 mm]
    \item
      Let $\cC$ be a small category of abelian groups that contains all $\Z/n\Z$ for $n\in\N_+$ and $\Q/\Z$
      as objects,
      let $\cI$ be the category of the  ordered set $\N_+$, (now ordered by $\mid$),
      let $*:\hat{S}\rightarrow\widehat{\str{S}}$ be an enlargement where $\hat{S}$ contains
      $\cC$ and $\cI$ as $\hat{S}$-small categories,
      and let $G:\cI\rightarrow\cC$ be the functor that sends $n\in\N_+$ to $\Z/n\Z$
      and a morphism $m\mid n$ in $\cI$ to $\Z/m\Z\xrightarrow{1\mapsto n/m}\Z/n\Z$.\\[4 mm]
      Then $\cI$ is filtered, and we can apply
      \ref{corinfinite}\ref{corfil} to $\varinjlim G$:
      The category $\str{\cI}$ is the category of the ordered set $\str{\N_+}$,
      and if $h\in\str{\N_+}\setminus\N$
      is any infinite natural number, then as $i_\infty$ we can take the number $(h!)$.
      The category $\str{\cC}$ is a category of abelian groups
      (with internal group-structure and with internal group homomorphisms as morphisms),
      and $\str{G}$ is the functor that sends a $k\in\str{\N_+}$ to $\str{\Z}/k\str{\Z}$.
      The corollary then tells us
      that to give a morphism from $\varinjlim G$ to an abelian group $A$ in $\cC$,
      we have to give an (internal) group homomorphism $\varphi$ from $\str{\Z}/h!\str{\Z}$ to $\str{A}$,
      such that for all $n\in\N_+$, the composition
      $\psi:\str{(\Z/n\Z)}\xrightarrow{1\mapsto h!/n}\str{\Z}/h!\str{\Z}\xrightarrow{\varphi}\str{A}$
      is of the form $\str{\varphi_n}$ for a $\varphi:\Z/n\Z\rightarrow A$.
      As this condition is obviously equivalent to the condition that $\psi(1)\in A$,
      we see that a morphism from $\varinjlim G$ to $A$ is given by a morphism
      $\varphi:\str{\Z}/h!\str{\Z}\rightarrow\str{A}$ such that for all $n\in\N_+$, we have
      $\varphi(h!/n)\in A\subseteq\str{A}$.\\[2 mm]
      Two such morphisms $\varphi_1$ and $\varphi_2$ correspond to the same morphism $\varinjlim G\rightarrow A$
      if and only if their compositions with all inclusions $\Z/n\Z\hookrightarrow\str{\Z}/h!\str{\Z}$
      coincide, i.e. if and only if $\varphi_1(h!/n)=\varphi_2(h!/n)\in A$ for all $n\in\N_+$.\\[2 mm]
      We have the isomorphism $\alpha:\varinjlim G\xrightarrow{\sim}\Q/\Z$,
      defined by the morphisms $\alpha_n:\Z/n\Z\xrightarrow{1\mapsto 1/n}\Q/\Z$ for $n\in\N_+$.
      To which morphisms $\varphi:\str{\Z}/h!\str{\Z}\rightarrow\str{\Q/\Z}$ does $\alpha$ correspond?
      --- By transfer, $\varphi$ is uniquely determined by $x:=\varphi(1)$,
      and by \eqref{eqpreshaveslim}, for every $n\in\N_+$, the following diagram in $\str{\cC}$ must commute:
      \[
        \xymatrix{
          {\str{\Z}/h!\str{\Z}} \ar[r]_{\varphi}^{1\mapsto x} & {\str{(\Q/\Z)}} \\
          {\Z/n\Z} \ar[u]^{1\mapsto h!/n} \ar[ru]_{\alpha_n} \\
        }
      \]
      This shows that we must have $1/n=\alpha_n(1)=(h!/n)\cdot x\in\str{(\Q/\Z)}$, i.e. $x=(1+mn)/h!$
      for an $m\in\str{\Z}$.
      Because this holds for every $n\in\N_+$, we get $x=(1+N)/h!$ with $N\in\str{\Z}$ a number divisible by all
      $n\in\N_+$.
  \end{enumerate}
\end{bsp}

\vspace{\abstand}

\section{Enlargements of additive and abelian categories}

\vspace{\abstand}

Next, we want to have a look at $\hat{S}$-small categories with more structure,
namely at \emph{additive} and \emph{abelian} categories.
We start by giving a formal definition:

\vspace{\abstand}

\begin{defi}
  Let $\hat{S}$ be a superstructure.
  \begin{enumerate}
    \item
      An \emph{additive $\hat{S}$-small category} is a pair
      $\langle\cA,P\rangle$,
      consisting of an $\hat{S}$-small category
      $\cA=\langle M,s,t,c\rangle$ and a set $P\subseteq M\times M\times M$,
      subject to the following conditions:
      \begin{enumerate}
        \item
          For all objects $A$, $B$ of $\cA$, the intersection
          $\bigl[\Mor{\cA}{A}{B}\times\Mor{\cA}{A}{B}\times\Mor{\cA}{A}{B}\bigr]\cap P$
          is a map $P_{A,B}:\Mor{\cA}{A}{B}\times\Mor{\cA}{A}{B}\rightarrow\Mor{\cA}{A}{B}$ that
          endows $\Mor{\cA}{A}{B}$ with the structure of an abelian group.
        \item
          For all objects $A$, $B$, $C$ of $\cA$, the map
          $\Mor{\cA}{B}{C}\times\Mor{\cA}{A}{B}\xrightarrow{\circ}\Mor{\cA}{A}{C}$
          given by $c$
          is bilinear with respect to the group-structures defined in (1).
        \item
          $\cA$ has a zero object.
        \item
          $\cA$ has arbitrary finite sums.
      \end{enumerate}
    \item
      An additive $\hat{S}$-small category $\langle\cA,P\rangle$ is called \emph{abelian}
      if it is abelian in the usual sense,
      i.e. if all morphisms in $\cA$ have a kernel and a cokernel
      and if for each morphism $f$, the canonical map
      $\text{coim}(f)\rightarrow\text{im}(f)$ is an isomorphism.
  \end{enumerate}
\end{defi}

\vspace{\abstand}

\begin{bem}
  \mbox{}\\
  \begin{enumerate}
    \item
      For a superstructure $\hat{S}$,
      an additive $\hat{S}$-small category is just an additive category in the usual sense whose
      underlying category is $\hat{S}$-small.
    \item
      Let $\cA$ be a small, additive category.
      Then there is a base set $S$ and an $\hat{S}$-small additive category which is isomorphic to $\cA$
      (compare \ref{bspssmall}\ref{bspsmallcat}).
  \end{enumerate}
\end{bem}

\vspace{\abstand}

\begin{bem}\label{bemex}
  Let $F:\cA\rightarrow\cB$ be an additive functor between abelian categories.
  Then $F$ is exact if and only if
  $F$ maps $\text{ker}(f)$ to $\text{ker}(Ff)$ for every morphism $f$ in $\cA$
  and if $F$ maps epimorphisms to epimorphisms:\\[2 mm]
  It is clear that the condition is necessary. If on the other hand
  $A\xrightarrow{f}B\xrightarrow{g}C$ is exact in $\cA$, then
  this is equivalent to the existence of a factorization
  \xymatrix{
    & A \ar[dr]^{f} \ar@{-->>}[d]_{\bar{f}} & & \\
    0 \ar[r] & {\text{ker}(g)} \ar[r] & B \ar[r]^{g} & {C} \\
  }\\
  with an epimorphism $\bar{f}$. Applying $F$ to this diagram and using $F\text{ker}(g)=\text{ker}(Fg)$
  and the fact that $F\bar{f}$ is an epimorphism then shows that
  $FA\xrightarrow{Ff}FB\xrightarrow{Fg}FC$ is exact in $\cB$,
  so the condition is also sufficient.
\end{bem}

\vspace{\abstand}

\begin{satz}\label{satzsternex}
  Let $*:\hat{S}\rightarrow\widehat{\str{S}}$ be an enlargement, and let $\langle\cA,P\rangle$
  be an additive $\hat{S}$-small category. Then:
  \begin{enumerate}
    \item\label{sternexadd}
      $\langle\str{\cA},\str{P}\rangle$ is an additive $\widehat{\str{S}}$-small category,
      and the functor $*$ is additive.
    \item\label{sternextrans}
      If $\varphi[X,Y,Z,f,g]$ is the formula
      $\bigl[
          \text{$X\xrightarrow{f}Y$ is a morphism in $\cA$ with kernel $Z\xrightarrow{g}X$}
        \bigr]
      $ in $\hat{S}$,
      then $\str{\varphi}[X,Y,Z,f,g]$ is the formula
      $\bigl[
          \text{$X\xrightarrow{f}Y$ is a morphism in $\str{\cA}$ with kernel $Z\xrightarrow{g}X$}
      \bigr]$ in $\widehat{\str{S}}$,
      and the analogous statement is true for ``cokernel'', ``image'' or ``coimage'' instead of ``kernel''.
    \item\label{sternexab}
      If in addition
      $\langle\cA,P\rangle$ is abelian, then so is $\langle\str{\cA},\str{P}\rangle$,
      and the functor $*$ is exact.
  \end{enumerate}
\end{satz}

\vspace{\abstand}

\begin{proof}
  By easy transfer, we see that $\langle\str{\cA},\str{P}\rangle$ satisfies
  conditions (1) and (2), and \ref{corfinlims} shows that it also satisfies
  conditions (3) and (4), so it is indeed an additive category.\\[1 mm]
  By \ref{correp}\ref{correpiii}, we know that $*$ maps the zero object of
  $\cA$ to the zero object of $\str{\cA}$
  and the sum $X\oplus Y$ of two objects $X,Y\in\cA$ to the sum of
  $\str{X}$ and $\str{Y}$ in $\str{\cA}$ --- this shows that $*$ is an \emph{additive} functor
  and therefore completes the proof of \ref{sternexadd}.\\[2 mm]
  To prove \ref{sternextrans}, we only have to formalize $\varphi$:
  \begin{multline*}
    \varphi[X,Y,Z,f,g]
    =\Bigl[
      \bigl(X,Y,Z\in\Ob(\cA)\wedge f\in\Mor{\cA}{X}{Y}\wedge g\in\Mor{\cA}{Z}{X}\bigr) \\
      \wedge\bigl(
        \forall A\in\Ob(\cA):
        \forall h\in\Mor{\cA}{A}{X}:
        [fh=0]
        \Rightarrow
        [\exists!h'\in\Mor{\cA}{A}{Z}:gh'=h]
      \bigr)
    \Bigr].
  \end{multline*}
  It is then clear (again by using \ref{corlemma}\ref{corobmor}) that $\str{\varphi}$ is what we claim in
  \ref{sternextrans}.\\[2 mm]
  Now let us assume that $\langle\cA,P\rangle$ is abelian.
  First, in $\str{\cA}$ kernels and cokernels of arbitrary morphisms exist because of \ref{corfinlims}.
  To prove that coimage and image of every morphism in $\str{\cA}$ are canonically isomorphic,
  we first write down the corresponding statement in $\cA$:
  \begin{multline*}
    \forall A,B,C,D\in\Ob(\cA):
    \forall f\in\Mor{\cA}{A}{B}: \\
    \forall\varphi\in\Mor{\cA}{A}{C}:
    \forall g\in\Mor{\cA}{C}{D}:
    \forall\psi\in\Mor{\cA}{D}{B}: \\
    \Bigl[
      \bigl(\text{$A\xrightarrow{\varphi}C$ is the coimage of $f$}\bigr)
      \wedge\bigl(\text{$D\xrightarrow{\psi}B$ is the image of $f$}\bigr)
      \wedge\bigl(\psi g\varphi=f\bigr)
    \Bigr] \\
    \Rightarrow
    \text{$g$ is an isomorphism}.
  \end{multline*}
  To see that the transfer of this statement is exactly what we want in $\str{\cA}$,
  we can apply \ref{sternextrans} for ``image'' and ``coimage''. This proves that $\str{\cA}$ is indeed
  an abelian category.\\[1mm]
  Finally, $*$ maps kernels to kernels because of \ref{correp}\ref{correpiii}
  and epimorphisms to epimorphisms because of \ref{trst}\ref{trepi},
  so \ref{bemex} shows that $*$ in this case is exact.
\end{proof}

\vspace{\abstand}

From now on, when talking about an additive or abelian $\hat{S}$-small category $\langle\cA,P\rangle$,
we will simply denote it by $\cA$ and drop $P$ from the notation.

\vspace{\abstand}

\begin{satz}\label{satzfunaddex}
  Let $*:\hat{S}\rightarrow\widehat{\str{S}}$ be an enlargement,
  let $\cA$ and $\cB$ be $\hat{S}$-small additive categories,
  and let $F:\cA\rightarrow\cB$ be a functor.
  \begin{enumerate}
    \item\label{funadd}
      If $F$ is additive, then $\str{F}:\str{\cA}\rightarrow\str{\cB}$ is also additive.
    \item\label{funex}
      If $\cA$ and $\cB$ are abelian
      and if $F$ is left exact (resp. right exact, resp. exact),
      then $\str{F}$ is also left exact (resp. right exact, resp. exact).
  \end{enumerate}
\end{satz}

\vspace{\abstand}

\begin{proof}
  To prove \ref{funadd}, we have to show that $\str{F}$ maps the zero object to the zero object and
  direct sums to direct sums:
  \[
    F(0_{\cA})=0_{\cB}
    \;\;\;\;\stackrel{\text{\tiny transfer}}{\Longrightarrow}\;\;\;\;
    (\str{F})(\underbrace{\str{0_{\cA}}}_{=0_{\str{\cA}}})=\underbrace{\str{0_{\cB}}}_{=0_{\str{\cB}}},
  \]
  and
  \begin{multline*}
    \forall X,Y\in\Ob(\cA):
    F(X\oplus Y)=F(X)\oplus F(Y) \\
    \stackrel{\ref{correp}\ref{correpii},\text{\tiny transfer}}{\Longrightarrow}\;\;\;\;
    \forall X,Y\in\Ob(\str{\cA}):
    (\str{F})(X\oplus Y)=(\str{F})(X)\oplus(\str{F})(Y).\\[1 mm]
  \end{multline*}
  Now let us assume that $\cA$ and $\cB$ are abelian and that $F$ is left exact.
  We have to show that $\str{F}$ then also is left exact, i.e. that $\str{F}$ maps kernels to kernels.
  But this follows, using \ref{corlemma}\ref{corobmor} and \ref{correp}\ref{correpiii}, by transfer of
  \[
    \forall A,B\in\Ob(\cA):
    \forall f\in\Mor{\cA}{A}{B}:
    F\text{ker}(f)=\text{ker}(Ff).
  \]
\end{proof}

\vspace{\abstand}

\begin{lemma}\label{lemmafunaddex}
  Let $*:\hat{S}\rightarrow\widehat{\str{S}}$ be an enlargement,
  and let $\cA$ and $\cB$ be $\hat{S}$-small categories.
  \begin{enumerate}
    \item\label{lemmafunadd}
      If $\cA$ and $\cB$ are additive, then the formula
      \[
        \varphi[F]\equiv
        \Bigl[
          \text{$F$ is an additive functor from $\cA$ to $\cB$}
        \Bigr]
      \]
      becomes
      \[
        \str{\varphi}[F]\equiv
        \Bigl[
          \text{$F$ is an additive functor from $\str{\cA}$ to $\str{\cB}$}
        \Bigr].
      \]
    \item\label{lemmafunex}
      If $\cA$ and $\cB$ are abelian, then the formula
      \[
        \varphi[F]\equiv
        \Bigl[
          \text{$F$ is an exact functor from $\cA$ to $\cB$}
        \Bigr]
      \]
      becomes
      \[
        \str{\varphi}[F]\equiv
        \Bigl[
          \text{$F$ is an exact functor from $\str{\cA}$ to $\str{\cB}$}
        \Bigr].
      \]
  \end{enumerate}
\end{lemma}

\vspace{\abstand}

\begin{proof}
  To prove \ref{lemmafunadd}, we can formalize $\varphi$ as follows:
  \[
    \varphi[F]
    =\Bigl[F\in\Ob(\cov{\cA}{\cB})\Bigr]
    \wedge\Bigl[F(0_{\cA})=0_{\cB}\Bigr]
    \wedge\Bigl[
      \forall A,B\in\Ob(\cA):
      F(A\oplus B)=F(A)\oplus F(B)
    \Bigr],
  \]
  and then the claim follows by transfer and using
  \ref{corlemma}\ref{corobmor} and \ref{correp}\ref{correpii}.
  To prove \ref{lemmafunex}, we use \ref{bemex} and write
  \begin{multline*}
    \varphi[F]
    =\Bigl[\text{$F$ is an additive functor from $\cA$ to $\cB$}\Bigr] \\
    \wedge\Bigl[\text{$F$ maps kernels to kernels and epimorphisms to epimorphisms}\Bigr],
  \end{multline*}
  and if we use \ref{lemmafunadd}, \ref{trst}\ref{trepi}, \ref{corlemma}\ref{corobmor} and \ref{correp}\ref{correpiii},
  the claim again follows by transfer.
\end{proof}

\vspace{\abstand}

Let $*:\hat{S}\rightarrow\widehat{\str{S}}$ be an enlargement,
let $R$ be a ring whose underlying set is an element of $\hat{S}$,
and let $\cA$ be an $\hat{S}$-small abelian category which is an abelian subcategory of the category of $R$-modules.

What does $\str{\cA}$ look like? --- Intuitively, we would think that $\str{\cA}$ is a subcategory of the category
of $\str{R}$-modules such that the objects of $\str{\cA}$ carry \emph{internal} $\str{R}$-module structures
and such that morphisms in $\str{\cA}$ are \emph{internal} $\str{R}$-linear maps. To prove this, we would like to reason
as follows: We take the statement $\bigl[\forall X\in\Ob(\cA):\text{$X$ is a set with an $R$-module structure}\bigr]$,
transfer it to get the desired result for objects of $\str{\cA}$
and then proceed similarly for morphisms.

Unfortunately, this reasoning is not completely correct, because it is not true that objects in $\cA$ ``are"
sets with an $R$-module structure. To see this, we must remember that the set of objects of $\cA$ was defined
as a certain set of \emph{morphisms} of $\cA$, namely the set of \emph{identity} morphisms.

To make matters more complicated, we do not even know what morphisms ``are": If $f$ is a morphism in $\cA$,
we do not really know what the \emph{set} $f$ is when we consider $f$ as an element of the superstructure $\hat{S}$,
because that depends on the way $\cA$ was made into an $\hat{S}$-small category.

Of course we could look at the exact definition of $\cA$
and then adjust the above reasoning to fit that exact description.
But we do not really want to do this. We do not want to have to know how exactly an $\hat{S}$-small category is
represented, because that would be rather tedious,
and in general there will be many ways to equip a given small category
with the structure of an $\hat{S}$-small category.

It seems better to use a \emph{categorical} approach instead of relying on the sets that represent the objects
of our category $\cA$. Instead of talking about the set that a given object ``is", we should rather look at the
\emph{functor} $G:\cA\rightarrow\text{Ens}$ that sends every object to its underlying set and then study the functor
$\str{G}$.

This is exactly what we are going to do in the next proposition (where we will make the additional assumption
that the functor $G$ is representable) to prove our initial intuition to be right.

\begin{satz}\label{satzrmodfun}
  Let $*:\hat{S}\rightarrow\widehat{\str{S}}$ be an enlargement,
  let $R$ be a ring whose underlying set is an element of $\hat{S}$,
  let $\rmod{R}$ be the category of $R$-modules,
  let $\cA$ be an $\hat{S}$-small abelian category,
  let $F:\cA\rightarrow\rmod{R}$ be an exact, fully faithful functor,
  and assume that there exists an object $\mR$ of $\cA$ such that
  there is an isomorphism $F\mR\xrightarrow{\psi}R$ of $R$-modules.
  Then we have:
  \begin{enumerate}
   \item\label{satzrmodfunii}
      Let $\rmod{R}(S)$ denote the full subcategory of $\rmod{R}$ consisting of $R$-modules whose underlying set
      is an element of $\hat{S}$.
      For every $A\in\Ob(\cA)$, the set $\Mor{\cA}{\mR}{A}$ carries
      an $R$-module structure induced by $F$ and $\psi$, and this defines an exact, fully faithful functor
      $G:=\Mor{\cA}{\mR}{\_}$ from $\cA$ to $\rmod{R}$ which is isomorphic to $F$ and which factorizes over
      $\rmod{R}(S)$.
    \item\label{satzrmodfuni}
      For every $A\in\Ob(\str{\cA})$, the set $\Mor{\str{\cA}}{\str{\mR}}{A}$
      carries an $\str{R}$-module structure induced by $\str{G}$,
      and this defines an exact, faithful functor
      $\tilde{G}:=\Mor{\str{\cA}}{\str{\mR}}{\_}$ from $\str{\cA}$ to $\rmod{\str{R}}$.
      Furthermore,
      for an object $M$ of $\rmod{R}(S)$, we have a canonical $\str{R}$-module structure on the set $\str{M}$
      and therefore get a functor $*:\rmod{R}(S)\rightarrow\rmod{\str{R}}$.
      Then the functors $(*\circ G)$ and $(\tilde{G}\circ *)$ from $\cA$ to $\rmod{\str{R}}$ are equal:
      \[
        \xymatrix@R=2mm@C=3mm{
          {\str{\cA}} \ar@{}[rrrrddd]|{=} \ar@{^{(}->}[rrrr]^{\tilde{G}} & & & & {\rmod{\str{R}}} & & & & & \\
          & & & & & & & & & \\
          &  & & & & & & & & \\
          {\cA} \ar@{^{(}->}@/_13mm/[rrrrrrrr]_{F} \ar@{^{(}->}[rrrr]_{G} \ar[uuu]^{*} & & & &
            {\rmod{R}(S)} \ar@{^{(}->}[rrrr] \ar[uuu]_{*} & & & & {\rmod{R}} \\
          & & \ar@{=>}^{\sim}[rrrr] & & & & & & & \\
        }
      \]
      This means that
      if we identify
      $\cA$ with an abelian subcategory of $\rmod{R}$ via $G$,
      we can identify $\str{\cA}$ with an abelian subcategory of $\rmod{\str{R}}$ by
      identifying an object $A$ of $\str{\cA}$ with the $\str{R}$-module $\Mor{\str{\cA}}{\str{\mR}}{A}$.
      Then for every $R$-module $M$ that is an object of $\cA$, $\str{M}$ is just the $\str{R}$-module $\str{M}$.
    \item\label{satzrmodfuniii}
      Let $\fin{\str{A}}$ denote the full subcategory of $\str{A}$ consisting of objects $A$ such that
      $\tilde{G}A$ is \emph{finite} (as a set).
      Then $\fin{\str{\cA}}$ is a Serre subcategory, i.e. it is closed with
      respect to subobjects, quotients and extensions.
  \end{enumerate}
\end{satz}

\vspace{\abstand}

\begin{proof}
  For an object $A$ of $\cA$,
  we define an $R$-module structure on $GA:=\Mor{\cA}{\mR}{A}$ by transport of structure via the bijections
  \[
      \Mor{\cA}{\mR}{A}\xrightarrow{F}
      \Mor{\rmod{R}}{F\mR}{FA}\xrightarrow{\psi^*}
      \Mor{\rmod{R}}{R}{FA} \xrightarrow{\text{\tiny{can}}}
      FA
  \]
  (of which the last two are also isomorphisms of $R$-modules).
  It is clear that $\psi^*$ induces an isomorphism between the functors $F$ and $G$,
  and because $F$ is exact and fully faithful, so is $G$.
  In addition to that, because $\cA$ is $S$-small, the set $GA$ is an element of $\hat{S}$,
  i.e. $G$ factorizes over $\rmod{R}(S)$. This proves \ref{satzrmodfunii}.\\[4 mm]
  Because of \ref{satzrmodfunii},
  for every $A\in\Ob(\cA)$, the set $GA$ carries an $R$-module structure,
  given by maps $a_A:GA\times GA\rightarrow GA$ and $s_A:R\times GA\rightarrow GA$.
  We thus get maps
  $a:\MorC{\cA}\times\MorC{\cA}\rightarrow\MorC{\cA}$ and
  $s:R\times\MorC{\cA}\rightarrow\MorC{\cA}$
  with the following properties:
  \begin{itemize}
    \item
      $\forall A\in\Ob(\cA):\forall x,y\in GA:a(x,y)\in GA$,
    \item
      $\forall A\in\Ob(\cA):\forall r\in R:\forall x\in GA:s(r,x)\in GA$,
    \item
      $\forall A\in\Ob(\cA):\text{$a$ and $s$ turn $GA$ into an $R$-module}$,
    \item
      \begin{minipage}[t]{14cm}
        \begin{multline*}
          \forall A,B\in\Ob(\cA):\forall\varphi\in\Mor{\cA}{A}{B}:
          \forall r\in R:
          \forall x,y\in GA: \\
          \varphi\circ\Bigl(a\bigl[s(r,x),y\bigr]\Bigr)
          =a\bigl[s(r,\varphi\circ x),\varphi\circ y\bigr]\in GB.
        \end{multline*}
      \end{minipage}
  \end{itemize}
  We therefore get maps
  $\str{a}:\MorC{\str{\cA}}\times\MorC{\str{\cA}}\rightarrow\MorC{\str{\cA}}$ and
  $\str{s}:\str{R}\times\MorC{\str{\cA}}\rightarrow\MorC{\str{\cA}}$
  which by transfer have the following properties
  (here we consider $G$ as a map from $\Ob(\cA)$ to $\MorC{\cA}$ which induces a map
  $\str{G}:\Ob(\str{\cA})\rightarrow\MorC{\str{\cA}}$):
  \begin{itemize}
    \item
      $\forall A\in\Ob(\str{\cA}):\str{G}(A)=\Mor{\str{\cA}}{\str{\mR}}{A}$,
    \item
      $\forall A\in\Ob(\cA):\forall x,y\in\str{G}(A):\str{a}(x,y)\in\str{G}(A)$,
    \item
      $\forall A\in\Ob(\cA):\forall r\in\str{R}:\forall x\in\str{G}(A):\str{s}(r,x)\in\str{G}(A)$,
    \item
      $\forall A\in\Ob(\cA):\text{$\str{a}$ and $\str{s}$ turn $\str{G}(A)$ into a $\str{R}$-module}$,
    \item
      \begin{minipage}[t]{14cm}
        \begin{multline*}
          \forall A,B\in\Ob(\str{\cA}):\forall\varphi\in\Mor{\str{\cA}}{A}{B}:
          \forall r\in\str{R}:
          \forall x,y\in\str{G}(A): \\
          \varphi\circ\Bigl(\str{a}\bigl[\str{s}(r,x),y\bigr]\Bigr)
          =\str{a}\bigl[\str{s}(r,\varphi\circ x),\varphi\circ y\bigr]\in\str{G}(B).
        \end{multline*}
      \end{minipage}
  \end{itemize}
  This shows that $\Mor{\str{\cA}}{\str{\mR}}{\_}$, together with $\str{a}$ and $\str{s}$,
  defines a functor $\tilde{G}:\str{\cA}\rightarrow\rmod{\str{R}}$.\\[3 mm]
  To see that $\tilde{G}$ is \emph{additive}, we note that because of the additivity of $G$ we have
  \[
    \forall A,B\in\Ob(\cA):\forall\varphi,\varphi'\in\Mor{\cA}{A}{B}:
    \forall x\in GA:
    (\varphi+\varphi')\circ x=s(\varphi\circ x,\varphi'\circ x),
  \]
  and this by transfer implies that $\tilde{G}$ is additive as well.\\[3 mm]
  Because $G$ is faithful, we have
  \[
    \forall A,B\in\Ob(\cA):
    \forall\varphi,\varphi'\in\Mor{\cA}{A}{B}:
    \Bigl(\bigl[\forall x\in GA:\varphi\circ x=\varphi'\circ x\in GB\bigr]
    \Rightarrow\varphi=\varphi'\Bigr),
  \]
  and transfer of this implies the faithfulness of $\tilde{G}$.\\[3 mm]
  To show that $\tilde{G}$ is \emph{exact},
  Note that for each exact sequence $A\rightarrow B\rightarrow C$ in $\cA$,
  the sequence $GA\rightarrow GB\rightarrow GC$ is exact in $\rmod{R}$.
  It follows by transfer from this that $\tilde{G}$ is exact.\\[3 mm]
  The last statement, namely that $\tilde{G}\circ *=*\circ G$,
  is immediately clear from the construction of $\tilde{G}$. This completes the proof of \ref{satzrmodfuni}.\\[4 mm]
  To prove \ref{satzrmodfuniii}, let $A\rightarrow B\rightarrow C$ be an exact sequence in $\str{\cA}$ such that
  $\tilde{G}A$ and $\tilde{G}C$ are finite. Because $\tilde{G}$ is exact,
  $\tilde{G}A\rightarrow\tilde{G}B\rightarrow\tilde{G}C$ is an exact sequence of $\str{R}$-modules,
  and it follows immediately that $\tilde{G}B$ is a finite set, i.e. $B$ belongs to $\fin{\str{\cA}}$ as well.
\end{proof}

\vspace{\abstand}

\begin{bem}\label{bemrmodfin}
  Proposition \ref{satzrmodfun} in particular applies to the case where $\cA$ is (equivalent to)
  the category of \emph{finitely generated} $R$-modules (for a ring $R$ whose underlying set is an element of
  $\hat{S}$, compare \ref{bspssmall}\ref{bsprmodfingen}). Also note that if it applies to a category $\cA$,
  then it also applies to any (full) exact subcategory $\cB\subseteq\cA$, even if $\cB$ does not contain the object $\mR$.
  For example, for $\cB$ we can take the category of all \emph{finite} $R$-modules.
  In particular, we can do this for the case $R=\Z$ and $\cB$ the category $\fin{\Abe}$ of finite abelian groups.
\end{bem}

\vspace{\abstand}

\def\ra{\rightarrow}

\section{Enlargements and derived functors}

\vspace{\abstand}

\begin{defi}
  A \emph{huge} set is a set of sets $N$ with the property that there is a universe $U$
  (see \ref{bem:universe}), containing an infinite set,
  such that the category $\Ens{N}$ is equivalent to the category $\Ens{U}$.
\end{defi}

\vspace{\abstand}

\begin{bem}
  Let $U$ be a universe.
  Then there is a base set $S$ and a set of sets $N$ in the superstructure $\hat{S}$
  such that the categories $\Ens{N}$ and $\Ens{U}$ are equivalent (compare \ref{bspssmall}\ref{bspsmallcat}).
  So in particular, $N$ is huge.
\end{bem}

\vspace{\abstand}

\begin{defi}
  Let $\hat{S}$ be a superstructure, and let $N\in\hat{S}$ be a huge set.
  Then by $\Ab{N}$ we denote the full subcategory of the category $\Abe$ of abelian groups
  consisting of groups $A$ whose underlying set is an element of $N$.\\[1 mm]
  Note that $\Ab{N}$ is an $\hat{S}$-small abelian category,
  and that
  if $\cA=\langle M,s,t,c\rangle$ is an $\hat{S}$-small abelian category
  with $M\subseteq N$, then $h_A=\Mor{\cA}{\_}{A}$ is an additive functor from
  $\cA$ to $\Ab{N}$ for every object $A$ in $\cA$.
\end{defi}

\vspace{\abstand}

\begin{satz}\label{satzab}
  Let $*:\hat{S}\rightarrow\widehat{\str{S}}$ be an enlargement,
  and let $N\in\hat{S}$ be a huge set. Then:
  \begin{enumerate}
    \item\label{satzabint}
      $\str{\Ab{N}}$ is the subcategory of $\Abe$
      whose objects are abelian groups in $\str{N}$ with an \emph{internal} group structure
      and whose morphisms are \emph{internal} group homomorphisms.
    \item\label{satzabex}
      The inclusion functor
      $\str{\Ab{N}}\hookrightarrow\Abe$ is \emph{exact}.
    \item\label{satzabinj}
      If $I$ is an injective object of $\str{\Ab{N}}$,
      then $I$ is also an injective object of $\Abe$.
    \end{enumerate}
\end{satz}

\vspace{\abstand}

\begin{proof}
  To prove \ref{satzabint}, we first compute the objects of $\str{\Ab{N}}$:\\
  \[
    \begin{array}{lcl}
      \Ob(\str{\Ab{N}}) &
      \stackrel{\ref{corlemma}\ref{corobmor}}{=} &
      \str{\Ob\Ab{N}} \\[2 mm]
      &= &
      \str{\bigl\{\langle A,m\rangle\bigl\vert
        A\in N,\
        \text{$m:A\times A\rightarrow A$ abelian group structure}
      \bigr\}} \\[2 mm]
      & = &
      \bigl\{\langle A,m\rangle\bigl\vert
        A\in\str{N},\
        \text{$m:A\times A\rightarrow A$ internal and abelian group structure}
      \bigr\}. \\
    \end{array}
  \]
  To compute the morphisms, we look at the following statement in $\hat{S}$:
  \[
    \forall A,B\in\Ob\Ab{N}:
    \Mor{\Ab{N}}{A}{B}
    =\bigl\{\varphi:A\rightarrow B\bigl\vert\text{$\varphi$ group homomorphism}\bigr\}.
  \]
  Transfer of this gives us (again using \ref{corlemma}\ref{corobmor})
  \[
    \forall A,B\in\Ob(\str{\Ab{N}}):
    \Mor{\str{\Ab{N}}}{A}{B}
    =\bigl\{\varphi:A\rightarrow B\bigl\vert\text{$\varphi$ internal group homomorphism}\bigr\}.
  \]
  Now we want to prove \ref{satzabex}. Denote the inclusion functor by $i$. First,
  it is easy to see that $i$ is additive. Next, look at the following statement in $\hat{S}$:
  \begin{multline*}
    \forall A,B,C\in\Ob\Ab{N}:
    \forall f\in\Mor{\Ab{N}}{A}{B}:
    \forall g\in\Mor{\Ab{N}}{B}{C}: \\
    \Bigl[
      \text{$A\xrightarrow{f}B\xrightarrow{g}C$ exact in $\Ab{N}$}
    \Bigr]
    \Leftrightarrow
    \Bigl[
      \bigl(
        gf=0
      \bigr)
      \wedge\bigl(
        \forall b\in B:
        [g(b)=0]\Rightarrow[\exists a\in A:f(a)=b]
      \bigr)
    \Bigr].
  \end{multline*}
  We can compute the transfer of this by using \ref{corlemma}\ref{corobmor} and \ref{lemmafunaddex}\ref{lemmafunex}:
  \begin{multline*}
    \forall A,B,C\in\Ob(\str{\Ab{N}}):
    \forall f\in\Mor{\str{\Ab{N}}}{A}{B}:
    \forall g\in\Mor{\str{\Ab{N}}}{B}{C}: \\
    \Bigl[
      \text{$A\xrightarrow{f}B\xrightarrow{g}C$ exact in $\str{\Ab{N}}$}
    \Bigr]
    \Leftrightarrow
    \underbrace{\Bigl[
      \bigl(
        gf=0
      \bigr)
      \wedge\bigl(
        \forall b\in B:
        [g(b)=0]\Rightarrow[\exists a\in A:f(a)=b]
      \bigr)
    \Bigr]}_{
    \Leftrightarrow\big[\text{$iA\xrightarrow{if}iB\xrightarrow{ig}iC$ exact in $\Abe$}\bigr]}.
  \end{multline*}
  This shows that the exactness of a sequence $A\xrightarrow{f}B\xrightarrow{g}C$ in $\str{\Ab{N}}$ is
  equivalent to the exactness of the sequence $iA\xrightarrow{if}iB\xrightarrow{ig}iC$ in $\Abe$
  and therefore proves \ref{satzabex}.\\[2 mm]
  To prove \ref{satzabinj} we use the well known fact (see for example \cite{tohoku}) that if
  $\cA$ is a full abelian subcategory of $\Abe$ that contains $\Z$,
  then $I$ is an injective object of $\cA$ if and only if $I$ is divisible. Because $N$ is huge,
  the category $\Ab{N}$ satisfies this condition,
  so the following statement is true in $\str{S}$:
  \[
    \forall I\in\Ob\Ab{N}:
    \forall n\in\N_+:
    \forall a\in I:
    \exists y\in I:
    x=n\cdot y.
  \]
  By transfer we get
  \[
    \forall I\in\Ob(\str{\Ab{N}}):
    \forall n\in\str{\N_+}:
    \forall a\in I:
    \exists y\in I:
    x=n\cdot y,
  \]
  which because of $\N_+\varsubsetneq\str{\N_+}$
  in particular implies that $I$ is divisible, i.e. an injective object of $\Abe$.
\end{proof}

\vspace{\abstand}

\begin{bem}\label{bemstrmod}
  Let $*:\hat{S}\rightarrow\widehat{\str{S}}$ be an enlargement,
  let $N\in\hat{S}$ be a huge set,
  let $R\in\hat{S}$ be a ring,
  let $\cA$ be the $\hat{S}$-small abelian category of $R$-modules $M$ with $M\in N$
  (resp. a subcategory of this category),
  and let $\cB$ be the abelian category of $\str{R}$-modules.
  Then it is easy to see (using arguments like in the proof of \ref{satzab})
  that $\str{\cA}$ is the subcategory of $\cB$ whose
  objects are $M\in\str{N}$ with \emph{internal} $\str{R}$-module structure
  and whose morphisms are \emph{internal} $\str{R}$-linear maps
  (resp. a subcategory of this category),
  and that the inclusion functor $\str{\cA}\hookrightarrow\cB$ is exact.
\end{bem}

\vspace{\abstand}

\begin{lemma}\label{lemmainjformula}
  Let $*:\hat{S}\rightarrow\widehat{\str{S}}$ be an enlargement,
  let $\cA$ be an $\hat{S}$-small abelian category,
  and consider the formula $\varphi[X]\equiv\Bigl[\text{$X$ is an injective object of $\cA$}\Bigr]$
  in $\hat{S}$.
  Then
  \[
    \str{\varphi[X]}=\Bigl[\text{$X$ is an injective object of $\str{\cA}$}\Bigr].
  \]
\end{lemma}

\vspace{\abstand}

\begin{proof}
  Let $\cA=\langle M,s,t,c\rangle$,
  and choose a huge $N\in\hat{S}$ containing $M$.
  Then
  \[
    \varphi[X]
    =\bigl(X\in\Ob(\cA)\bigr)\wedge
    \bigl(\text{$h_X:\cA\rightarrow\Ab{N}$ is exact}\bigr).
  \]
  We compute $\str{\varphi}[X]$ with \ref{lemmafunaddex}\ref{lemmafunex}:
  \[
    \str{\varphi}[X]
    =\bigl(X\in\Ob(\str{\cA})\bigr)\wedge
    \bigl(\text{$\Mor{\str{\cA}}{\_}{X}:\str{\cA}\rightarrow\str{\Ab{N}}$ is exact}\bigr).
  \]
  Now because of \ref{satzab}\ref{satzabex}, we know that the exactness of $\Mor{\str{\cA}}{\_}{X}$
  is the same as the exactness of $h_X:\str{\cA}\rightarrow\Abe$,
  and by definition this exactness is equivalent to $X$ being injective, so the lemma follows.
\end{proof}

\vspace{\abstand}

\begin{cor}\label{corinj}
  Let $*:\hat{S}\rightarrow\widehat{\str{S}}$ be an enlargement,
  and let $\cA$ be an $\hat{S}$-small abelian category.
  Then $*$ maps injective (resp. projective) objects of $\cA$ to injective (resp. projective)
  objects of $\str{\cA}$,
  and if $\cA$ has enough injectives (resp. projectives), then so has $\str{\cA}$.
\end{cor}

\vspace{\abstand}

\begin{proof}
  The corollary follows easily from \ref{lemmainjformula}: If $I$ is an injective object of $\cA$, then
  $\varphi[I]$ is true, so that by transfer $\str{\varphi}[\str{I}]$ is also true,
  i.e. $\str{I}$ is an injective object of $\str{\cA}$.
  And if $\cA$ has enough injectives, then the following statement is true:
  \[
    \forall A\in\Ob(\cA):
    \exists I\in\Ob(\cA):
    \exists f\in\Mor{\cA}{A}{I}:
    \Bigl[\text{$f$ monomorphism}\Bigr]\wedge\varphi[I].
  \]
  Then by transfer and because of \ref{trst}\ref{trmono} and \ref{corlemma}\ref{corobmor},
  it follows that also $\str{\cA}$ has enough injectives.
  --- The proof for projectives is analogous.
\end{proof}

\vspace{\abstand}

\begin{cor}\label{corfuninjinj}
  Let $*:\hat{S}\rightarrow\widehat{\str{S}}$ be an enlargement,
  and let $\cA$ and $\cB$ be $\hat{S}$-small abelian category,
  and let $F:\cA\rightarrow\cB$ be a functor that maps injectives to injectives
  (resp. projectives to projectives).
  Then $\str{F}:\str{\cA}\rightarrow\str{\cB}$ also maps injectives to injectives
  (resp. projectives to projectives).
\end{cor}

\vspace{\abstand}

\begin{proof}
  For injectives, this follows immediately from \ref{lemmainjformula} by transfer of the statement
  \[
    \forall A\in\Ob(\cA):
    \bigl[\text{$A$ injective}\bigr]
    \Rightarrow
    \bigl[\text{$FA$ injective}\bigr],
  \]
  and for projectives, it is analogous.
\end{proof}

\vspace{\abstand}

\begin{bem}
  Let $\cA$ be a Grothendieck category with a generator
  $U\in\cA$. For an object $A\in\cA$ we define the cardinality
  of $A$ by $|A|:=|\{A'\subset A\}|$. Let $\kappa$ be an inaccessible
  (i.e. regular and limit) cardinal number with $\kappa>|U|$
  and $\kappa>|Hom(U,U)|$ and let $\cA^{<\kappa}$ be the full
  subcategory of $\cA$ of objects $A\in\cA$ with $|A|<\kappa$.
  Further let $S$ be a set with $|S|>\kappa$.

  \begin{flushright}\begin{minipage}{15cm}

    \vspace{\abstand}

    \begin{nonumberlemma}\mbox{}\\[-4mm]
      \begin{enumerate}
      \item The category $\cA^{<\kappa}$ is a small abelian category
        with enough injectives objects and the inclusion functor
        $i:\cA^{<\kappa}\ra\cA$ is exact and respects injective
        objects.
      \item The category $\cA^{<\kappa}$ is an $\hat{S}$-small category
      \end{enumerate}
    \end{nonumberlemma}

    \vspace{\abstand}

    \begin{proof}[Sketch of proof]
      Each object of $\cA^{<\kappa}$ is a quotient of the object $\sum_{\kappa}U$. Therefore $\cA^{<\kappa}$ is small.
      If $0\ra A' \ra A \ra A'' \ra 0$ is a short exact sequence in $\cA$, than the inequalities $|A'|\leq |A|$ and
      $|A''|\leq |A|$ hold. We further have $|A\times B|\leq 2^{\alpha^{max(|A|,|B|)}}$. So the category $\cA{<\kappa}$
      is a abelian subcategory of $\cA$ and the inclusion functor is exact. The concrete construction of injective
      resolutions in \cite{tohoku} and the choice of $\kappa$ show that $\cA{<\kappa}$ has enough injectives. Because
      $U$ and all subobjects of $U$ are in  $\cA{<\kappa}$ the inclusion functor $i$ respects injectives (comp.
      \cite[Lemma 1 in 1.10 ]{tohoku}).
    \end{proof}

  \end{minipage}\end{flushright}

\end{bem}

\vspace{\abstand}

\begin{cor}\label{corrflf}
  Let $*:\hat{S}\rightarrow\widehat{\str{S}}$ be an enlargement,
  and let $F:\cA\rightarrow\cB$ be an additive functor between $\hat{S}$-small abelian categories.
  \begin{enumerate}
    \item\label{corrf}
      If $\cA$ has enough injectives and if $F$ is left exact, then for all $i\geq 0$,
      the right derived functors
      $\RF{F}{i}:\cA\rightarrow\cB$ and $\RF{(\str{F})}{i}:\str{\cA}\rightarrow\str{\cB}$
      exist, and the following diagram of functors commutes:
      \[
        \xymatrix@C=2cm{
          {\cA} \ar[r]^{*} \ar[d]_{\RF{F}{i}} &
          {\str{\cA}} \ar@{=}[r] \ar[d]^{\str{(\RF{F}{i})}} &
          {\str{\cA}} \ar[d]^{\RF{(\str{F})}{i}} \\
          {\cB} \ar[r]_{*} & {\str{\cB}} \ar@{=}[r] & {\str{\cB}} \\
        }
      \]
    \item\label{corlf}
      If $\cA$ has enough projectives and if $F$ is right exact, then for all $i\geq 0$
      the left derived functors
      $\LF{F}{i}:\cA\rightarrow\cB$ and $\LF{(\str{F})}{i}:\str{\cA}\rightarrow\str{\cB}$
      exist, and the following diagram of functors commutes:
      \[
        \xymatrix@C=2cm{
          {\cA} \ar[r]^{*} \ar[d]_{\LF{F}{i}} &
          {\str{\cA}} \ar@{=}[r] \ar[d]^{\str{(\LF{F}{i})}} &
          {\str{\cA}} \ar[d]^{\LF{(\str{F})}{i}} \\
          {\cB} \ar[r]_{*} & {\str{\cB}} \ar@{=}[r] & {\str{\cB}} \\
        }
      \]
  \end{enumerate}
\end{cor}

\vspace{\abstand}

\begin{proof}
  We only give a proof for \ref{corrf}, the proof for \ref{corlf} is analogous.
  It is well known that the right derived functors of a left exact functor exist if the source category has enough
  injectives. Now $\cA$ has enough injectives by assumption, so the $\RF{F}{i}$ exist,
  and $\str{F}$ is left exact by \ref{satzfunaddex}\ref{funex} and $\str{\cA}$ has enough injectives
  because of \ref{corinj}, so the $\RF{(\str{F})}{i}$ also exist.\\[1 mm]
  The left square commutes because of \ref{satzfunctors}\ref{sfunctorsiii}.
  To show that the right square commutes, look at the following
  statement in $\hat{S}$:
  \begin{multline}\label{eqinj}
    \forall A\in\Ob(\cA):
    \forall\text{$(A\xrightarrow{\varepsilon}I^0\xrightarrow{\delta^0}I^1\xrightarrow{\delta^1}\ldots
      \xrightarrow{\delta^i}I^{i+1})$ injective resolution of $A$ in $\cA$}:\\
      \RF{F}{i}A=\Koh{I^{i-1}\xrightarrow{\delta^{i-1}}I^i\xrightarrow{\delta^i}I^{i+1}}.
  \end{multline}
  What is the transfer of this statement? To answer this question,
  we first note that the transfer of
  \[
    \varphi[X,Y,Z,f,g]\equiv\Bigl[\text{$X\xrightarrow{f}Y\xrightarrow{g}Z$ is an exact sequence in $\cA$}\Bigr]
  \]
  is
  \[
    \str{\varphi}[X,Y,Z,f,g]=\Bigl[\text{$X\xrightarrow{f}Y\xrightarrow{g}Z$ is an exact sequence in $\str{\cA}$}\Bigr].
  \]
  (To see this, we write the exactness in terms of equalities between kernels and cokernels etc. and then use
  \ref{correp}\ref{correpii}.)
  From \ref{lemmainjformula} we know that the transfer of the statement
  $\varphi[X]\equiv\Bigl[\text{$X$ is injective in $\cA$}\Bigr]$
  is $\str{\varphi}[X]=\Bigl[\text{$X$ is injective in $\str{\cA}$}\Bigr]$,
  and combining these two results with \ref{satzsternex}\ref{sternextrans} we get as the transfer of \eqref{eqinj}:
  \begin{multline*}
    \forall A\in\Ob(\str{\cA}):
    \forall\text{$(A\xrightarrow{\varepsilon}I^0\xrightarrow{\delta^0}I^1\xrightarrow{\delta^1}\ldots
      \xrightarrow{\delta^i}I^{i+1})$ injective resolution of $A$ in $\str{\cA}$}:\\
      \str{(\RF{F}{i})}A=\underbrace{\Koh{I^{i-1}\xrightarrow{\delta^{i-1}}I^i\xrightarrow{\delta^i}I^{i+1}}}_
      {=\RF{(\str{F})}{i}A}.
  \end{multline*}
  This proves that the right square also commutes.
\end{proof}

\vspace{\abstand}

\newcommand\Kom[1]{{\text{Kom}\left(#1\right)}}
\newcommand\Komb[1]{{\text{Kom}^b \left(#1\right)}}

\def\cT{\mathcal T}
\def\cK{\mathcal K}
\def\cF{\mathcal F}
\def\cE{\mathcal E}
\def\cZ{\mathcal Z}
\def\ra{\rightarrow}

\section{Enlargement of triangulated and derived categories}

\vspace{\abstand}

In this section we want to investigate what happens if we enlarge triangulated
and derived categories. From now on, the proofs will be a little bit shorter,
and we will not give all the details.

\vspace{\abstand}

\noindent
Let $\cT$ be an additive category with an automorphism $\Sigma: \cT\ra\cT$.
\begin{itemize}
  \item A triangle in $\cT$ is a triple $(f,g,h)$ of morphisms in
    $\cT$ of the form $A\xrightarrow{f} B\xrightarrow{g} C\xrightarrow{h}\Sigma A$
    with $A,B,C\in\cT$.
  \item The category $\cT$ together with the automorphism $\Sigma$ and a class of triangles $\Delta$
    is called a triangulated category if it fulfils certain axioms (compare for
    example \cite{Weib}).
\end{itemize}

\vspace{\abstand}

Now let $\hat{S}$ be a superstructure, $(\cT,\Sigma,\Delta)$ a triangulated
category with $\cT$ $\hat{S}$--small, and $*:\hat{S}\ra\widehat{\str{S}}$ be
an enlargement. In particular this means that the set $\Delta$ is also
$\hat{S}$--small.

\noindent
We get
\begin{satz}\label{triangstar}
  The triple $(\str{\cT}, \str{\Sigma},\str{\Delta})$ is again
  a triangulated category.
\end{satz}

\vspace{\abstand}

\begin{proof}
  Follows easily from the transfer principle.
\end{proof}

\vspace{\abstand}

\noindent
We also get
\begin{satz}\label{trianfuncstar}
  If $F:\cT\ra\cT'$ is a functor of triangulated categories then
  $\str{F}:\str{\cT}\ra\str{\cT'}$ is also a functor of triangulated categories.
\end{satz}

\vspace{\abstand}

\begin{proof}
  Follows again easily from the transfer principle.
\end{proof}

\vspace{\abstand}

\noindent
Now let us look at triangulated subcategories.

\begin{defi}
  A triangulated subcategory $\cD$ of a triangulated category $\cT$ is called \emph{thick}
  iff for all objects $x,y\in\Ob(\cT)$, we have the implication $x\oplus y\in\Ob(\cD)\Rightarrow x\in\Ob(\cD)$.
\end{defi}

\vspace{\abstand}

\noindent
We have already seen that $*$ is exact. The next statement reflects this on a triangulated level.
\begin{satz}\label{starthick}
  Let $\cT$ be an $\hat{S}$--small triangulated category and $\cD\subset\cT$ be a
  thick subcategory. Then $\str{\cD}\subset\str{\cT}$ is again thick
  and we have a canonical isomorphism
  $$ \str{\cT}/\str{\cD} \xrightarrow{\sim} \str{(\cT/\cD)}$$
  of triangulated categories.
\end{satz}

\vspace{\abstand}

\begin{proof}
  That $\str{\cD}\subset\str{\cT}$ is again a thick subcategory follows from the
  definition, the transfer principle, and corollary \ref{correp} (ii). The canonical functor $\cT\ra\cT/\cD$
  induces a functor $\str{\cT}\ra\str{(\cT/\cD)}$ and by transfer the subcategory
  $\str{\cD}$ is the full subcategory of objects which are mapped to zero. So we
  get a functor $\str{\cT}/\str{\cD} \ra \str{(\cT/\cD)}$. Now we have the
  following general and well known fact about triangulated categories:
  \begin{lemma}
    Let $\cT$ and $\cT'$ be triangulated categories and $F:\cT\ra\cT'$ be a functor of
    triangulated categories with
    \begin{enumerate}
    \item $F$ is essentially surjective on objects
    \item for all $A',B'\in\Ob(\cT')$ and $f\in\Mor{\cT'}{A'}{B'}$ there are
      $A,B\in\Ob(\cT)$, $f\in\Mor{\cT}{A}{B}$ and isomorphism $g_A:F(A)\ra A'$
      and $g_B:F(B)\ra B'$ such that the diagram
   $$\xymatrix{
      F(A)\ar[r]^{F(f)} \ar[d]^{g_A} & F(B) \ar[d]^{g_B} \\
      A' \ar[r]^f                      & B }$$
     is commutative.
   \end{enumerate}
   Then the canonical functor $\cT/ker(F) \ra \cT'$ is an isomorphism of triangulated categories.
  \end{lemma}
  By the transfer principle and \ref{satzstarcat} both properties are fulfilled for the functor
  $\str{\cT}\ra \str{(\cT/\cD)}$ and so the proposition follows.
\end{proof}

\vspace{\abstand}

Now we want to consider the derived category of an abelian category.
The derived category can be obtained in three steps. So let $\cA$ be an
abelian category.  First
we look at the category $\Kom{\cA}$ of complexes in the abelian category. Then
we identify morphisms of complexes which are homotopic and get
the homotopy--category of complexes $\cK(\cA)$. This is a triangulated category
and there we divide out the thick subcategory of all complexes which are
quasi--isomorphic to the zero complex to get the derived
category $\cD(\cA)$. Now we look step by step what happens under enlargements.
For that let $\cA$ be an $\hat{S}$--small abelian category. Let
$\cZ$ be the category which has as objects $\Z$ and for each
$i\leq j$ exactly one morphism from $i$ to $j$. For a functor
from $\cZ$ to a category we denote by $d_i$ the image of the
unique morphism from $i$ to $i+1$. The category $\Kom{\cA}$ of
complexes in $\cA$ is the category of functors from $\cZ$ to
$\cA$ with the property that for all $i\in\cZ$  the equality
$d_{j+1}\circ d_{j}=0$ holds. Therefore $\str{\Kom{\cA}}$ is the category of internal
functors from $\str{\cZ}$ to $\str{\cA}$ with $d_{j+1}\circ d_{j}=0$ for all $i\in\str{\cZ}$.
We call the objects of this category \emph{$*$--complexes}.
The category $\str{\cZ}$ has as objects $\str{\Z}$ and again
for each $i\leq j$ exactly one morphism from $i$ to $j$. Because
$\Z$ lies inside of $\str{\Z}$ we get a restriction  functor from
$\str{\Kom{\cA}}$ to $\Kom{\str{\cA}}$. By the transfer principle
we get $\str{\cK(\cA})$ out of $\str{\Kom{\cA}}$ if we identify
morphisms which are internal homotopic. Because internal homotopies
in $\str{\Kom{\cA}}$ become homotopies in $\Kom{\str{\cA}}$ this functor
induces a functor from $\str{\cK(\cA)}$ to $\cK(\str{\cA})$.
With \ref{starthick} one can see easily that we get a functor
$$res_{\cA}:\str{\cD(\cA)}\ra\cD(\str{\cA}).$$

For an abelian category $\cA$ we denote by $\cD^b(\cA)$ the full triangulated subcategory
of $\cD(\cA)$ of all complexes $A^{\bullet}\in\cD(\cA)$ such that there is an $i\in\N$
with $\forall j\in \Z: |j|>i \Rightarrow h^i(A^{\bullet})\simeq 0$. The full subcategory
of $\Kom{\cA}$ of all complexes $A^{\bullet}\in\Kom{\cA}$ with $A^j\simeq 0$ for all
$j\in\Z$ with $|j|>i$ for a specific $i\in\N$ is denoted by $\Komb{\cA}$.
\vspace{\abstand}

\begin{bem}
  It is not possible to define $res_{\cA}$ for the bounded derived category, i.e. a functor
$\str{\cD^b(\cA)}\ra\cD^b(\str{\cA})$, because the *--complexes
in $\str{\cD^b(\cA)}$ are not bounded by a standard natural number.
\end{bem}

\vspace{\abstand}

\noindent
By transfer we get for each $i\in\str{\Z}$ a functor
$$\bullet^i:\str{\Kom\cA}\ra \str{\cA}$$
and
$$\str{h}^i: \str{\cD(\cA)} \ra \str{\cA}.$$

\vspace{\abstand}

\noindent
We summarize this in the following
\begin{satz}\label{inculdestar}
  Let $\cA$ be an $\hat{S}$--small abelian category.
  Then for each $i\in\Z$, we have the following commutative diagram:
  \[
    \xymatrix{
      &  \str{\cA} \ar[rr]^{\str{i}} \ar@{=}'[d][dd]   &    &
        \str{\cD(\cA)}  \ar@{->}'[d] ^{res_{\cA}}[dd]  \ar[rr]^{\str{h^i}} &
      &   \str{\cA} \ar@{=}[dd] \\
      \cA \ar[ur]^{*} \ar[rd]_{*} \ar[rr]^>>>>>>>i &  &
        \cD(\cA) \ar[ur]^{*} \ar[dr]^{R*} \ar[rr]^>>>>>>>{h^i}  &   &
        \cA \ar[ru]^{*} \ar[rd]^{*}  &  \\
      & \str{\cA}  \ar[rr]_i &   & \cD(\str{\cA}) \ar[rr]^{h^i}   &  &  \str{\cA}
    }
  \]
  Here $R*$ denotes the right derived functor of the exact functor $*:\cA\rightarrow\str{\cA}$.
\end{satz}

\vspace{\abstand}

\begin{proof}
  This follows immediately from the constructions.
\end{proof}

\vspace{\abstand}

In general we do not have a non-trivial functor from $\cD(\str{\cA})$ to $\str{\cD(\cA)}$. But we do have
one if we restrict ourselves to bounded complexes. Then we even get more. Namely we can identify $\cD^b(\str{\cA})$
with a full but external subcategory of $\str{\cD(\cA)}$.

\begin{satz}\label{derfib}
  Let $\cA$ be an $\hat{S}$--small abelian category.
  There is a canonical fully faithful functor
  $$ \cD^b(\str{\cA}) \ra \str{\cD(\cA)}$$
  which identifies $\cD^b(\str{\cA})$ with the (external!) triangulated subcategory
  $\str{\cD^{fb}(\cA)}$ of all *--complexes $A\in\str{\cD(\cA)}$ such that
  there is an $i\in\N$ with
  $$\forall j\in\str{\Z}:|j|>i\Rightarrow h^i(A)\simeq 0 \text{ in } \str{\cA}.$$
  The functor $res_{\cA}:\str{\cD(\cA)}\ra\cD(\str{\cA})$ restricts to a functor
  $\str{\cD^{fb}(\cA)}\ra\cD^b(\str{\cA})$ which is an inverse of the above functor.
\end{satz}

\vspace{\abstand}

\begin{proof}
  An object in $\cD^b(\str{\cA})$ is represented by a complex $A^{\bullet}\in\Kom{\str{\cA}}$ with
  $A^j\simeq 0$ for all $j\in\Z$ with $|j|>i$ for a certain $i\in\Z$.
  Because there are only finitely many $j\in\Z$ with $A^j\not= 0$,
  we can continue $A^{\bullet}$ to an internal functor from $\str{\cZ}$ to $\str{\cA}$ (compare \ref{satzinternal}).
  The same is true for homotopies and quasiisomorphisms. So the statement follows from this.
\end{proof}

\vspace{\abstand}

\begin{satz}\label{starder}
  Let $\cA$ and $\cB$ be $\hat{S}$--small abelian categories, and we assume that $\cA$ has enough
  injective objects. Further let $F:\cA\ra\cB$ be a left exact functor. Then the functors in the
  following diagram exist and make it commutative.
  $$\xymatrix{
     \cD^{+}(\cA) \ar[r]^{RF} \ar[d]_{R*}  &  \cD^{+}(\cB) \ar[d]_{R*}  \\
     \cD^{+}(\str{\cA}) \ar[r]^{R\str{F}}  &  \cD^{+}(\str{\cB})
     }$$
\end{satz}

\vspace{\abstand}

\begin{proof}
  This can be proven in the same way as proposition \ref{corrflf}.
\end{proof}

\begin{satz}\label{derfibsatz}
   Let $\cA$ and $\cB$ be $\hat{S}$--small abelian categories, and we assume that $\cA$ has enough
  injective objects. Further let $F:\cA\ra\cB$ be a left exact functor of finite cohomological dimension.
  The functor $\str{RF}$ restricts to a functor
  $$ \str{\cD^{fb}(\cA)}\ra\str{\cD^{fb}(\cB)}$$
  and we get a commutative diagram
  $$\xymatrix{
    \str{\cD^{fb}(\cA)}\ar[rr]^{\str{RF}_{|\str{\cD^{fb}(\cA)}}}\ar[d]^{\wr}_{res_{\cA}}    &   &\str{ \cD^{fb}(\cB) \ar[d]^{\wr}_{res_{\cB}}}  \\
      \cD^b(\str{\cA}) \ar[rr]^{R\str{F}}                                        &    &  \cD^b(\str{\cB})
             }$$
\end{satz}

\vspace{\abstand}

\begin{proof}
  We denote by $\str{\text{Kom}^{fb}(\cA)}$ the category of all $*$--complexes $A^{\bullet}\in\str{\Kom{\cA}}$
  such that
  there is an $i\in\N$ with
  $$\forall j\in\str{\Z}:|j|>i\Rightarrow A^{i}\simeq 0 \text{ in } \str{\cA}.$$
  By transfer we have for each $*$--complex $A^{\bullet}\in\str{\text{Kom}^{fb}(\cA)}$
  a resolution in $\str{\text{Kom}^{fb}(\cA)}$ $r:A^{\bullet}\ra I^{\bullet}_{A^{\bullet}}$
  with $I^{\bullet}_{A^{\bullet}}\in\str{\text{Kom}^{fb}(\cA)}$ such that
  for all $ i\in\N$ the object $I^i_{A^{\bullet}}$ is $\str{F}$--acyclic
  and $\str{RF}=\str{F}(I^{\bullet}_{A^{\bullet}})$. We see that $res_{\cA}(I^{\bullet}_{A^{\bullet}})$
  has also acyclic components and so
  $R\str{F}(res_{\cA}(I^{\bullet}_{A^{\bullet}}))=\str{F}(res_{\cA}(I^{\bullet}_{A^{\bullet}})$
  and the proposition follows.
\end{proof}

\vspace{\abstand}

\begin{bem}
  It is easy to see that if $\cA$ and $\cB$ are $\hat{S}$-small abelian categories,
  if $\cA$ has enough \emph{projective} objects,
  and if $F:\cA\rightarrow\cB$ is a \emph{right} exact functor (of finite cohomological dimension),
  the obvious statements which are analogous to \ref{starder} and \ref{derfibsatz} hold.
\end{bem}

\vspace{\abstand}

\def\cF{\mathcal F}
\def\cE{\mathcal E}
\def\ra{\rightarrow}
\newcommand\Hom[3]{\mathrm{Hom}_{#1}(#2,#3)}

\section{Enlargements of fibred categories}

\vspace{\abstand}

For each ring in the superstructure $\hat{S}$ we have the
$\hat{S}$--small category of finitely generated modules over the
ring. This gives us also
for each standard ring (i.e. of the form $\str{R}$ for a ring
$R$ in $\hat{S}$) an internal category of modules in $\widehat{\str{S}}$.
But we also want to have a category of modules for an internal
but nonstandard ring.  Furthermore we would like to have
functorial properties with respect to morphisms between the rings.
For this we show in this section
that the notion of fibred categories behaves well
under enlargements. How this solves the above problem we will see
in example \ref{ringbeispiel}.

\vspace{\abstand}

For the theory of fibred categories we refer to \cite[exposé 6]{SGAI} and for
additive, abelian, triangulated and derived fibrations to \cite[exposé 17]{SGA4III}.

\vspace{\abstand}

\noindent
First we recall the notations and definitions
of fibred categories relevant for us.

\begin{defi}
  \mbox{}\\[-2mm]
  \begin{enumerate}
    \item
      Let $p:\cF\ra\cE$ be a functor,
      and let $S$ be an object in $\cE$.  The category $\cF_S$ with
      $$Ob(\cF_S):=\{X\in\cF|p(X)=S\}\text{ and }
           \Mor{\cF_S}{X}{Y}:=\{f\in\Mor{\cF}{X}{Y}|p(f)=\id{S}\}$$
      is called the \emph{fibre of $\cF$ in $S$}.
    \item
      Let $p:\cF\ra\cE$ be a functor,
      and let $\alpha:X\ra Y$ be a morphism in $\cF$.
      Then $\alpha$ is called \emph{cartesian}
      if for all $X'\in \cF$ and for all morphisms $u:X'\ra Y$
      such that there is a factorisation $p(u)=p(\alpha)\circ \beta$
      there is a unique $\overline{u}\in \Mor{\cF}{X'}{X}$ with 
      $u=\alpha\circ \overline{u}$ and $p(\overline{u})=\beta$.
    \item
      A functor $p:\cF\ra\cE$ is called
      a \emph{fibration} if for every morphism $\alpha:S'\ra S$ in
      $\cE$ and for every $Y\in\cF_{S}$,
      there is a cartesian morphism $f:X\ra Y$ with $p(f)=\alpha$.
      The category $\cF$ is also called a category \emph{fibred over $\cE$}.
  \end{enumerate}
\end{defi}

\vspace{\abstand}

\begin{bem}
  The definition of a cartesian morphism differs slighly from
  \cite[VI.5.1]{SGAI} but compare \cite[VI.6.11]{SGAI}.
\end{bem}

\vspace{\abstand}

\noindent
Now let $*:\hat{S}\rightarrow\widehat{\str{S}}$ be an enlargement.

\begin{satz}\label{enlfib}
  Let $p:\cF\ra\cE$ be a fibration of $\hat{S}$--small
  categories. Then $\str{p}:\str{\cF}\ra\str{\cE}$ is also
  a fibration and for all $S\in\cE$ we have
  $\str{(\cF_S)}=\str{\cF}_{\str{S}}$.
\end{satz}

\vspace{\abstand}

\begin{proof}
  The last statement follows immediately from the definition.
  Now let $\cF=\langle M,s_M,t_M,c_M\rangle$ and
  $\cE=\langle N,s_N,t_N,c_N\rangle$.  Then the functor is just
  a map $p:M\ra N$.  For a morphism in $\cF$ to be
  cartesian we have the formula
  \begin{multline*}
    \varphi[f]\equiv f\in M \wedge \forall X'\in\cF\;
    \forall u\in\Mor{\cF}{X'}{t_M(f)}:\exists \beta\in
    N:p(u)=p(f)\circ \beta \Rightarrow\\
    \exists ! \overline{u}\in\Mor{\cF}{X'}{s_M(f)}:
    u=f\circ\overline{u} \land p(\overline{u})=\beta
  \end{multline*}
  Now we have
\begin{multline*}
    \str\varphi[f]\equiv f\in \str M \wedge \forall X'\in\str\cF\;
    \forall u\in\Mor{\str\cF}{X'}{\str t_M(f)}:\exists \beta\in
    \str N:\str p(u)=\str p(f)\circ \beta \Rightarrow\\
    \exists ! \overline{u}\in\Mor{\str\cF}{X'}{\str s_M(f)}:
    u=f\circ\overline{u} \land \str p(\overline{u})=\beta
  \end{multline*}
  But this means just that $f$ is a cartesian morphism in $\str{\cF}$.
  Now the theorem follows by applying the transfer principle.
\end{proof}

\vspace{\abstand}

If all fibres of a fibration have a certain formal property,
then the fibres of the enlargement of the fibration have the same property.
We don't want to work this out in detail but we just want to
give one example.

\begin{satz}
  Let $p:\cF\ra\cE$ be a fibration of $\hat{S}$--small
  categories and we assume that for all $S\in\cE$ the category
  $\cF_S$ has fibre products. Then we also have for
  all $S\in\str{\cE}$ fibre products in $\str{\cF}_S$.
\end{satz}

\vspace{\abstand}

\begin{proof}
  This can be proven in the same way as in corollary \ref{corfinlim}
\end{proof}

\vspace{\abstand}

Now we want to add some more structures to the category $\cF$.
For that we recall some definitions from \cite{SGA4III}.
For a functor $p:\cA\ra\cB$, a morphism $f$ in $\cB$ and objects
$X,Y\in\cA$ we define $\Hom{f}{X}{Y}:=\{g\in\Hom{\cA}{X}{Y}|p(g)=f\}.$

\vspace{\abstand}

\begin{defi}
  An \emph{additive category over a category $\cB$} is a category $\cA$
  with a functor $p:\cA\ra\cB$ and with abelian group structures
  on all $\Hom{f}{X}{Y}$ for all morphisms $f$ in $\cB$ and all
  objects $X,Y\in\cA$ such that:
  \begin{enumerate}
  \item The composition in $\cA$ is biadditive.
  \item The fibres of the functor become additive categories.
  \end{enumerate}
\end{defi}

\vspace{\abstand}

\begin{defi}
  An \emph{abelian category over a category $\cB$} is an additive
  category $F:\cA\ra\cB$ such that the fibres are abelian
  categories and such that for all $f:x\ra y$ in $\cB$,
  the bifunctor $$\Hom{f}{X}{Y}:\cA_x^{op}\times\cA_y\ra\mathrm{Ab}$$
  is left exact in both $X$ and $Y$.
\end{defi}

\vspace{\abstand}

\begin{defi}
  A \emph{triangulated category over a category $\cB$} is an additive category
  $F:\cA\ra\cB$ over $\cB$ together with a functor $T:\cA\ra\cA$ and
  a class of triangles $\Delta_x$ in $\cA_x$ for all $x\in\cB$ such that:
  \begin{enumerate}
  \item T is an additive $B$--automorphism of $\cA$.
  \item For all $x\in\cB$ the triple $(\cA_x, T_x, \Delta_x)$ is a triangulated
    category.
  \item For all $x,y\in\cB$ and triangles $(X,Y,Z,u,v,w)\in\Delta_x,
    (X',Y',Z',u',v',w')\in\Delta_y$ and all commutative diagrams
    $$\xymatrix{
      X \ar[r]^u \ar[d]^f & Y \ar[r]^v \ar[d]^j & Z \ar[r]^w   & T_x X
               \ar[d]^f\\
      X' \ar[r]^{u'}      & Y' \ar[r]^{v'}      & {Z'}\ar[r]^{w'} & T_y X' }
    $$
    there is an $h:Z\ra Z'$ such that this will become a morphism of triangles.
  \end{enumerate}
\end{defi}

\vspace{\abstand}

\noindent
It is straightforward to enlarge such categories.
We omit a precise statement and just state the following

\begin{satz}\label{abfib}
  The notion of small additive (abelian/triangulated) category over another small category
  behaves well under enlargement.
\end{satz}

\vspace{\abstand}

\noindent
Moreover we get
\begin{satz}\label{injfib}
  Let $p:\cA\ra\cB$ be an abelian fibration of $\hat{S}$--small categories such
  that for all $B\in\cB$ the abelian category $\cA_B$ has enough injective/projective
  objects. Then for all $B\in\str{\cB}$ the abelian category $\str{\cA}_B$ has
  enough injective/projective objects.
\end{satz}

\vspace{\abstand}

\begin{proof}
  This can be proved in the same way as corollary \ref{corinj}.
\end{proof}

\vspace{\abstand}

\begin{bsp}\label{ringbeispiel}
  Let $\cE$ be the category of finite rings in $\hat{S}$ from example \ref{bsplim}.
  Now let $\cF$ be the category of pairs $(R,M)$ where $R$ is a finite Ring and
  $M$ is a finitely generated $R$--module. A morphism from an object
  $(R,M)$ to an object $(R',M')$ is by definition a pair $(\varphi, f)$
  where $\varphi$ is just a ring homomorphism $\varphi:R\ra R'$ and
  $f$ is just an $R$--linear module homomorphism from $M$ to $M'_{R}$.
  Here $M'_{R}$ denotes the module which is as abelian group the same as
  $M'$ and gets its $R$--module structure via $\varphi:R\ra R'$.
  This gives a category with the obvious composition. We have the
  natural functor $\cF\ra\cE$ which maps the pair $(R,M)$ to $R$ and
  a morphism $(\varphi,f)$ to $\varphi$. One can easily check that this gives
  an abelian fibred category over the category of finite rings.
  Now let $\str{\cF}\ra\str{\cE}$ be the corresponding
  fibration in $\widehat{\str{S}}$.  If $R$ is an object in $\str{\cE}$ we call
  the fibre $\str{\cF_R}$ the category of *--finitely generated $R$--modules.
  The category $\str{\cE}$ can easily be identified with the category
  of *--finite internal rings in $\hat{S}$. Furthermore $\str{\cF}_R$ can
  be identified with a category of specific internal $R$--modules in $\widehat{\str{S}}$.
\end{bsp}

\vspace{\abstand}

\begin{bem}\label{bemzweicat}
  A motivation for the previous definitions is the following proposition from \cite[XVII: Prop. 2.2.13]{SGA4III}

  \begin{flushright}\begin{minipage}{15cm}

      \vspace{\abstand}

      \begin{nonumbersatz}
    There is an equivalence between
    \begin{enumerate}
    \item the 2-category of pseudo--functors from a category $\cB^{op}$
      to the category of additive (abelian/triangulated) categories with additive
      (left exact/triangulated) functors
    \item the 2--category of additive (abelian/triangulated) fibred categories
      over $\cB$.
    \end{enumerate}
  \end{nonumbersatz}
  \end{minipage}\end{flushright}

  \vspace{\abstand}

  The notion of a 2--category also behaves well under enlargement.
  Now let $(Cat)$ be a $\hat{S}$--small 2--category of small categories in $\hat{S}$.
  Then one can get by the transfer principle from the previous proposition the following

  \begin{flushright}\begin{minipage}{15cm}

      \vspace{\abstand}

      \begin{nonumbersatz}
            Let $\cB$ be a $\hat{S}$--small category. There is an equivalence between
    \begin{enumerate}
    \item the 2-category of internal pseudo--functors from a category $\str{\cB^{op}}$
      to the category of additive (abelian/triangulated) categories in $\str{(Cat)}$ with additive
      (left exact/triangulated) functors.
    \item the 2--category of additive (abelian/triangulated) fibred internal categories
      over $\str{\cB}$.
    \end{enumerate}
  \end{nonumbersatz}
  \end{minipage}\end{flushright}
\end{bem}

\vspace{\abstand}

Finally we want to have a fibred version of proposition \ref{derfib} and \ref{derfibsatz}.
For details on derived fibrations of an abelian fibration and cofibrations we refer again to
\cite[expose 17]{SGA4III}.

\vspace{\abstand}

\begin{satz}\label{fibversion}
  Let $\cA\ra\cB$ be an abelian cofibration which is bounded derivable,
  and furthermore we assume that there is an $n_0\in\N$ such that for
  all morphisms $f:X\ra Y$ in $\cB$ the cohomological dimension of
  the left exact functor $f_*:\cA_X\ra\cA_Y$ is less or equal to $n_0$.
  Then $\str{\cA}\ra\str{\cB}$ is also an abelian cofibration which is
  bounded derivable and for all $f:X\ra Y$ in $\str{\cB}$ the cohomological
  dimension of the left exact functor $f_*:\str{\cA}_X\ra\str{\cA}_Y$
  is again less or equal to $n_0$. Further we define $\str{\cD^{fb}}(\cA/\cB)$
  as the full subcategory of $\str{\cD}(\cA/\cB)$ which consists of those *complexes which
  are bounded by a standard natural number. Then $\str{\cD^{fb}}(\cA/\cB)$
  defines also a triangulated cofibration over $\str{\cB}$ and we have
  an canonical isomorphism of triangulated fibrations
  $$\str{\cD^{fb}}(\cA/\cB) \xrightarrow{\sim} \cD^b(\str{\cA}/\str{\cB}).$$
\end{satz}

\vspace{\abstand}

\begin{proof}
  By the transfer principle it follows that $\str{\cA}/\str{\cB}$ is again derivable
  and the statement about the cohomological dimension. The last property
  shows that $\str{\cD^{fb}}(\cA/\cB)$ is also cofibred over $\str{\cB}$.
  The last statement of the proposition can be shown in a similar way as
  in \ref{derfib} and \ref{derfibsatz}.
\end{proof}

\vspace{\abstand}

\def\app{\angle}

\appendix
\section{Superstructures and Enlargements}

\vspace{\abstand}

In this appendix we give a short introduction to enlargements and list (without proofs) some of their basic properties.
Details can be found in \cite{lanrog} and \cite{loeb}.

\vspace{\abstand}

\noindent
For a set $M$, let $\cP(M)$ denote the power set of $M$, i.e. the set of all subsets of $M$.

\vspace{\abstand}

\begin{defi}(Superstructure)\\
  Let $S$ be an infinite set whose elements are no sets.
  Such a set we call a \emph{base set}, and its elements we call \emph{base elements}.
  We define $\hat{S}$, the \emph{superstructure over $S$}, as follows:
  \[
    \hat{S}:=\bigcup_{n=0}^\infty S_n
    \;\;\;\;\;\;
    \text{where $S_0:=S$ and $\forall n\geq 1:S_n:=S_{n-1}\cup\cP(S_{n-1})$.}
  \]
\end{defi}

\vspace{\abstand}

\begin{bem}\label{bem:universe}
  Recall from \cite[I.0]{SGA4I} that a \emph{universe} is a set $U$ satisfying the following conditions:
  \begin{enumerate}
    \item
      $y\in x\in U\Rightarrow y\in U$,
    \item
      $x,y\in U\Rightarrow\{x,y\}\in U$,
    \item
      if $x\in U$ is a set, then $\cP(x)\in U$,
    \item
      if $(x_i)_{i\in I\in U}$ is a family with $x_i\in U$ for all $i$, then
      $\bigcup_{i\in I}x_i\in U$.
  \end{enumerate}
  Every superstructure satisfies axioms (i), (ii) and (iii) of a universe, but not (iv):
  Let $\hat{S}$ be a superstructure.
  \begin{enumerate}
    \item
      Let $y\in x\in\hat{S}$. Because $x$ is a set and the base elements are no sets,
      there exists a uniquely determined $n\in\N_+$ such that $x\in S_n\setminus S_{n-1}$, i.e.
      $x\in\cP(S_{n-1})$. So $x$ is a subset of $S_{n-1}$ and $y$ an element of $S_{n-1}\subseteq\hat{S}$.
    \item
      Let $x,y\in\hat{S}$. Then there is a $n\in\N_0$ such that $x,y\in S_n$, i.e.
      $\{x,y\}\in\cP(S_n)\subseteq S_{n+1}\subseteq\hat{S}$.
    \item
      Let $x\in\hat{S}$ be a set. Then as in (i), $x\subseteq S_n$
      for a suitable $n\in\N_0$,
      and it follows $\cP(x)\subseteq\cP(S_n)\subseteq S_{n+1}$,
      i.e. $\cP(x)\in S_{n+2}\subseteq\hat{S}$.
    \item
      To see that (iv) does no hold, let $I$ be a countable infinite subset of $S$,
      choose a bijection $\varphi:I\rightarrow\N_0$ and put $x_i:=S_{\varphi(i)}\in\hat{S}$.
      (Note that for all $n$, the set $S_n$ is an element of $\cP(S_n)\subseteq S_{n+1}\subseteq\hat{S}$.)
      Then $\bigcup_{i\in I}x_i=\hat{S}$, which is \emph{not} an element of $\hat{S}$.
  \end{enumerate}
\end{bem}

\vspace{\abstand}

Note that even though superstructures fail to be universes, axiom (iv) holds under the additional assumption
that all sets $x_i$ of the family are contained in a fixed set $U\in\hat{S}$, and in praxis all families considered will
be of this sort. In addition to that, we have the following

\begin{satz}\label{satzUsmallSsmall}
  Let $U$ be a universe,
  and let $\cC$ be an $U$-small category.
  Then there exists a base set $S$ such that $\cC$ is equivalent to an $\hat{S}$-small category.
\end{satz}

\vspace{\abstand}

\begin{proof}
  Let $M\in U$ denote the set of morphisms of $\cC$. Choose a base set $S$ whose cardinality is greater than that of $M$,
  and choose an injection $M\hookrightarrow S$. Then the construction given after \ref{deficat} obviously produces
  an $\hat{S}$-small category that is isomorphic to $\cC$.
\end{proof}

\vspace{\abstand}

In the superstructure $\hat{S}$ to a base set $S$, we will find most of the mathematical objects of interest
related to $S$:
First of all, for sets $M,N\in\hat{S}$, the product set $M\times N$ is again an element of $\hat{S}$ when
we identify an ordered pair $\langle a,b\rangle$ for $a\in M$, $b\in N$ with the set
$\{a,\{a,b\}\}$,
and for sets $M_1,\ldots,M_n\in\hat{S}$, the product set
$M_1\times\ldots\times M_n:=(M_1\times\ldots\times M_{n-1})\times M_n$
is an element of $\hat{S}$.
Therefore, relations between two sets $M,N\in\hat{S}$ and in particular functions from $M$ to $N$
are again elements of $\hat{S}$.

For example, if $S$ contains the set of real numbers $\R$, then $\hat{S}$ will contain the sets $\R^n$ for
$n\in\N_+$ as well as functions between subsets of $\R^n$ and $\R^m$, the set of continuous functions between such sets
or the set of differentiable functions and so on.

\vspace{\abstand}

\begin{defi}
  Let $\hat{S}$ be a superstructure, and let $f,x\in\hat{S}$. Then
  we define the element $(f\app x)$ of $\hat{S}$ as follows:
  if $f$ is a set that contains $\langle x,y\rangle$ for exactly
  one $y\in\hat{S}$, then $(f\app x):=y$, otherwise $(f\app x):=\emptyset$.\\[2 mm]
  (So if $f:A\rightarrow B$ is a function with $A,B\in\hat{S}$ and if $x\in A$,
  then $(f\app x)=f(x)$.)
\end{defi}

\vspace{\abstand}

\noindent
Now we want to define what \emph{terms}, \emph{formulas} and \emph{statements} in a given superstructure are:

\begin{defi}
  Let $\hat{S}$ be a superstructure.
  \begin{enumerate}
    \item
      A \emph{term in $\hat{S}$} is the following
      \begin{enumerate}
        \item
          An element of $\hat{S}$ is a term, a so called \emph{constant}
          (note that $\hat{S}$ itself is not an element of $\hat{S}$ and therefore not a constant).
        \item
          A \emph{variable}, i.e. a symbol that is not an element of $\hat{S}$,
          is a term; usually we will denote variables by
          $x,y,x_1,x_2,X,Y$ and so on.
        \item
          If $s$ and $t$ are terms, then $\langle s,t\rangle$ and $(s\app t)$ are terms.
        \item
          A sequence of symbols is a term if and only if it can be built by
          (1), (2) and (3) in finitely many steps.
      \end{enumerate}
      Note that terms that do not contain variables are elements of $\hat{S}$.
    \item
      A \emph{formula in $\hat{S}$} is the following:
      \begin{enumerate}
        \item
          If $s$ and $t$ are terms, then $s=t$ and $s\in t$ are formulas.
        \item
          If $\psi$ is a formula, then $\neg\psi$ is a formula.
        \item
          If $\psi$ and $\chi$ are formulas, then $(\psi\wedge\chi)$ is a formula.
        \item
          If $x$ is a variable, $t$ is a term in which $x$ does not appear and $\psi$ is a formula,
          then $(\forall x\in t)\psi$ is a formula.
        \item
          A sequence of symbols is a formula if and only if it can be built by (1), (2), (3) and (4)
          in finitely many steps.
      \end{enumerate}
    \item
      Let $\varphi$ be a formula in $\hat{S}$,
      and let $X$ be a variable that occurs at position $j$ of
      $\varphi$ (remember that a formula is a sequence of symbols).
      We say that this occurrence of $X$ is \emph{bound} if position $j$ is part of a subsequence of
      $\varphi$ of the form
      $(\forall X\in t)\psi$
      for a term $t$ and a formula $\psi$. Otherwise, we call the
      occurrence of $X$ \emph{free}.
    \item
      Let $\varphi$ be a formula in $\hat{S}$.
      We say that a variable $X$ is a \emph{free variable of
      $\varphi$}, if there is a free occurrence of $X$ at a
      position of $\varphi$. If $X_1,\ldots,X_n$ are the free
      variables of $\varphi$, we often write
      $\varphi[X_1,\ldots,X_n]$ instead of $\varphi$.
    \item
      If $\varphi[X_1,\ldots,X_n]$ is a formula in $\hat{S}$ with free variables $X_1,\ldots,X_n$
      and if $\tau_1,\ldots,\tau_n\in\hat{S}$ are constants,
      then $\varphi[\tau_1,\ldots,\tau_n]$ is the statement in $\hat{S}$
      that we get when we replace any \emph{free} occurrence
      of $X_i$ in $\varphi$ by $\tau_i$.
    \item
      A \emph{statement in $\hat{S}$} is a formula that has no
      free variables.
  \end{enumerate}
\end{defi}

\vspace{\abstand}

\begin{satz}\label{satzformula}
  Let $\varphi$ be a statement in a superstructure $\hat{S}$.
  \begin{enumerate}
    \item
      If $\varphi$ contains none of the logical symbols $\neg$, $\wedge$ or $\forall$,
      then $\varphi$ is either of the form $s=t$ or of the form $s\in t$ with terms
      $s,t\in\hat{S}$ that contain no variables and are therefore elements of $\hat{S}$.
    \item
      If $\varphi$ contains at least one of the symbols $\neg$, $\wedge$ or $\forall$,
      then it is of one and only one of the following forms:
      \begin{enumerate}
        \item
          $\varphi=\neg\psi$ with a statement $\psi$,
        \item
          $\varphi=(\psi_1\wedge\psi_2)$ with statements $\psi_1$ and $\psi_2$.
        \item
          $\varphi=(\forall X\in t)\psi$ with a term $t$ without variables and
          a formula $\psi$ which can only contain $X$ as a free variable.
      \end{enumerate}
  \end{enumerate}
\end{satz}

\vspace{\abstand}

\begin{defi}
  Let $\varphi$ be a statement in a superstructure $\hat{S}$.
  We will call $\varphi$
  \emph{true} or \emph{valid} under the following conditions,
  depending on the structure of $\varphi$ in the sense of \ref{satzformula}:
  \begin{enumerate}
    \item
      If $\varphi$ is of the form $s=t$ or $s\in t$ with terms $s$ and $t$ that are elements of $\hat{S}$,
      then $\varphi$ is \emph{true} iff $s$ equals $t$ in the first case and iff $s$ is an element of $t$ in
      the second case.
    \item
      If $\varphi=\neg\psi$ for a statement $\psi$, then $\varphi$ is \emph{true} iff
      $\psi$ is not true.
    \item
      If $\varphi=(\psi_1\wedge\psi_2)$ for statements $\psi_1$ and $\psi_2$,
      then $\varphi$ is \emph{true} iff $\psi_1$ and $\psi_2$ both are true.
    \item
      If $\varphi=(\forall X\in t)\psi$ for a term $t$ which is an element of $\hat{S}$ and a formula
      $\psi$ which can only contain $X$ as a free variable, then we distinguish between the following three cases:
      \begin{enumerate}
        \item
          If $t$ is a set and if $X$ is a free variable of $\psi$, then
          $\varphi$ is \emph{true} iff the statement $\psi[\tau]$ is true for all $\tau\in t$.
        \item
          If $t$ is a set and if $\psi$ is a statement, then
          $\varphi$ is \emph{true} iff $\psi$ is true.
        \item
          If $t$ is no set, then $\varphi$ is \emph{true}.
      \end{enumerate}
  \end{enumerate}
\end{defi}

\vspace{\abstand}

\begin{defi}\label{defienlargement}(Enlargement) \\
  Let $*:\hat{S}\rightarrow\hat{W}$ be a map between superstructures.
  For $\tau\in\hat{S}$ we denote the image of $\tau$ under $*$ by $\str{\tau}$,
  and for a formula $\varphi$ in $\hat{S}$,
  we define $\str{\varphi}$ to be the formula in $\hat{W}$ that we get when we replace any constant
  $\tau$ occurring in $\varphi$ by $\str{\tau}$.\\[1mm]
  We call $*$ an \emph{enlargement} if the following conditions hold:
  \begin{enumerate}
    \item
      $\str{S}=W$.
      (Because of this property, we will often write $*:\hat{S}\rightarrow\widehat{\str{S}}$.)
    \item
      $\forall s\in S:\str{s}=s.$
    \item(\emph{transfer principle})\\
      If $\varphi$ is a statement in $\hat{S}$, then $\varphi$ is true
      iff $\str{\varphi}$ is true.
    \item(\emph{saturation principle})\label{saturation}\\
      Put $\cI:=\bigcup_{A\in\hat{S}\setminus S}\str{A}\subseteq\hat{W}$,
      let $I$ be a nonempty set whose cardinality is not bigger than that of $\hat{S}$,
      and let $\{U_i\}_{i\in I}$ be a family of nonempty sets $U_i\in\cI$
      with the property that for all \emph{finite} subsets $ J\subseteq I$,
      the intersection $\bigcap_{j\in J}U_j$ is nonempty. Then $\bigcap_{i\in I}U_i\neq\emptyset$.
  \end{enumerate}
  If $*$ is an enlargement, we call the elements of $\cI$ the \emph{internal elements} of $\hat{W}$.
\end{defi}

\vspace{\abstand}

\begin{thm}\label{thmenlargement}
  For any base set $S$, there exists an enlargement $*:\hat{S}\rightarrow\widehat{\str{S}}$.
\end{thm}

\vspace{\abstand}

An easy corollary of the transfer principle is that $*$ is injective: If $\sigma$ and $\tau$
are elements of $\hat{S}$ and $\varphi$ denotes the statement $\bigl[\sigma=\tau\bigr]$, then
$\str{\varphi}\equiv\bigl[\str{\sigma}=\str{\tau}\bigr]$, i.e. $\str{\sigma}=\str{\tau}$ implies
$\sigma=\tau$ by transfer.

From now on, in building formulas, we will
as well use the symbols $\vee$, $\Rightarrow$, $\Leftrightarrow$, $\exists$ and $\exists!$
as obvious abbreviations for the corresponding constructions involving only the symbols $\neg$, $\wedge$
and $\forall$. We will also often write $f(x)$ instead of $(f\app x)$ and $\forall X\in t:\psi$
instead of $(\forall X\in t)\psi$ to make formulas more readable.

\vspace{\abstand}

\begin{satz}\label{satzapp}
  Let $*:\hat{S}\rightarrow\widehat{\str{S}}$ be an enlargement.
  \begin{enumerate}
    \item
      If $M\in\hat{S}$ is a set, then $\str{M}$ is also a set.
    \item\label{satzappfin}
      If $A\in\hat{S}$ is a set, then
      $\{\str{a}\vert a\in A\}\subseteq\str{A}$ with equality if and only if $A$ is \emph{finite}.
      In particular, if $A\subset S$ is finite, then $\str{A}=A$.
    \item
      $\str\emptyset=\emptyset$.
    \item
      If $A,B\in\hat{S}$ are sets, then
      \[
        \str{(A\cup B)}=\str{A}\cup\str{B},\
        \str{(A\cap B)}=\str{A}\cap\str{B},\
        \str{(A\setminus B)}=\str{A}\setminus\str{B},\
        \str{(A\times B)}=\str{A}\times\str{B}.
      \]
    \item
      If $A_1,\ldots,A_n\in\hat{S}$ are sets and $a_1\in A_1,\ldots,a_n\in A_n$, then
      \[
        \str{\langle a_1,\ldots,a_n\rangle}
        =\langle\str{a_1},\ldots,\str{a_n}\rangle.
      \]
    \item
      Let $A\in\hat{S}$ be a set and $\varphi[X_1,\ldots,X_n]$ be a formula in $\hat{S}$. Then\\[1mm]
      \mbox{}\hspace{3mm}$
        \str{\bigl\{\langle a_1,\ldots,a_n\rangle\in A^n\bigl\vert\text{$\varphi[a_1,\ldots,a_n]$ is true}\bigr\}}
        =\bigl\{\langle a_1,\ldots,a_n\rangle\in(\str{A})^n\bigl\vert
        \text{$\str{\varphi}[a_1,\ldots,a_n]$ is true}\bigr\}
      $.
    \item
      If $A\in\hat{S}$ is a set, then
      \[
        \str{[\cP(A)]}
        =\bigl\{B\in\cP(\str{A})\bigl\vert\text{$B$ internal}\bigr\}.
      \]
    \item\label{satzappviii}
      For $f,x\in\hat{S}$, we have $\str{(f\app x)}=(\str{f}\app\str{x})$,
      and if
      $A,B\in\hat{S}$ are sets and $f:A\rightarrow B$ is a function, then
      $\str{f}$ is a function from $\str{A}$ to $\str{B}$.
    \item
      If $A,B,C\in\hat{S}$ are sets and $g:A\rightarrow B$ and $f:B\rightarrow C$ are maps,
      then
      $\str(f\circ g)=\str{f}\circ\str{g}$.
    \item
      If $A,B\in\hat{S}$, then
      \[
        \str{\left[B^A\right]}
        =\bigl\{f:\str{A}\rightarrow\str{B}\bigl\vert\text{$f$ internal}\bigr\}
      \]
      (here $B^A$ denotes as usual the set of maps from $A$ to $B$).
    \item
      If $A\in\hat{S}$ is an ordered set (resp. partial ordered set, totally ordered set, group, abelian group,
      ring, commutative ring, commutative ring with unit, field, totally ordered field),
      then $\str{A}$ also is an ordered set (resp. partial ordered set, \ldots).
  \end{enumerate}
\end{satz}

\vspace{\abstand}

\begin{satz}\label{satzinternal}
  Let $*:\hat{S}\rightarrow\widehat{\str{S}}$ be an enlargement.
  \begin{enumerate}
    \item
      Let $A,B\in\widehat{\str{S}}$ be \emph{internal} sets, i.e. sets which are internal elements of $\widehat{\str{S}}$.
      Then $A\cup B$, $A\cap B$, $A\setminus B$ and $A\times B$
      are also internal sets.
    \item
      Let $B\in\widehat{\str{S}}$ be an \emph{internal} set,
      and let $\varphi[X_1,\ldots,X_n]$ be an \emph{internal} formula in $\widehat{\str{S}}$,
      i.e. a formula in which all constants are internal elements of $\widehat{\str{S}}$. Then the set
      \[
        \bigl\{\langle b_1,\ldots,b_n\rangle\in B^n\bigl\vert\text{$\varphi[b_1,\ldots,b_n]$ is true}\bigr\}
      \]
      is internal.
    \item
      If $a\in A\in\widehat{\str{S}}$ and $A$ is an internal set, then $a$ is internal.
    \item\label{satzinternalfinite}
      If $A\in\hat{S}\setminus S$ is a \emph{finite} set,
      and $B\in\widehat{\str{S}}$ is an \emph{internal} set,
      then any map from $\str{A}=A$ to $B$ is internal.
  \end{enumerate}

\end{satz}

\vspace{\abstand}



\end{document}